\documentclass[12pt,reqno,tbtags]{amsart}
\usepackage{color}
\newcommand{\revised}[1]{\textcolor{black}{#1}}
\usepackage{amssymb}
\usepackage{mathrsfs,citesort}
\def\mathcal{\mathscr}
\usepackage[dvipdfm]{graphicx}

\nopagebreak[3]

\numberwithin{equation}{section}

\let\Sum=\sum
\DeclareMathOperator*{\textsum}{{\textstyle\Sum}}
\def\sum{\textsum}
\let\Prod=\prod
\DeclareMathOperator*{\textprod}{{\textstyle\Prod}}
\def\prod{\textprod}
\let\Bigcap=\bigcap
\DeclareMathOperator*{\textbigcap}{{\textstyle\Bigcap}}
\def\bigcap{\textbigcap}

\newtheorem{theorem}{Theorem}[section]
\newtheorem{proposition}{Proposition}[section]
\newtheorem{lemma}{Lemma}[section]
\newtheorem{corollary}{Corollary}[section]
   
\theoremstyle{definition}
   \newtheorem{definition}{Definition}[section]
   
\theoremstyle{remark}
   \newtheorem*{remark}{\rm \textit{Remark}}

\makeatletter
\newcommand{\thmskip}{\vspace{.5\baselineskip\@plus.2\baselineskip
                                    \@minus.2\baselineskip}\par\noindent}
\makeatother

\makeatletter
\def\@startSsection#1#2#3#4#5#6{%
 \if@noskipsec \leavevmode \fi
 \par \@tempskipa #4\relax
 \@afterindenttrue
 \ifdim \@tempskipa <\z@ \@tempskipa -\@tempskipa \@afterindentfalse\fi
 \if@nobreak \everypar{}\else
     \addpenalty\@secpenalty\addvspace\@tempskipa\fi
 \@ifstar{\@dblarg{\@sect{#1}{\@m}{#3}{#4}{#5}{#6}}}%
        {\@dblarg{\@sect{#1}{#2}{#3}{#4}{#5}{\noindent\S#6}}}%
}

\def\@secnumfont{}
\def\section{\@startSsection{section}{1}%
  \z@{.7\linespacing\@plus\linespacing}{.5\linespacing}%
  {\normalfont\bfseries}}
\makeatother

\newcommand{\QED}{~$\quad\square$}
\newcommand{\nablaA}{\nabla_{\!\!\A}}
\newcommand{\DeltaA}{\Delta_{\!\A}}
\newcommand{\A}{\textbf{\textit{A}}}
\newcommand{\B}{\textbf{\textit{B}}}
\newcommand{\J}{\textbf{\textit{J}}}
\newcommand{\M}{\mathrm{M}}

\newcommand{\LG}{\textrm{L}}
\newcommand{\CG}{\textrm{C}}
\newcommand{\R}{\textbf{\textit{R}}}
\newcommand{\C}{\textbf{\textit{C}}}

\newcommand{\embedding}{\hookrightarrow}
\newcommand{\Proof}{\par\noindent {\it Proof}. }
\newcommand{\Proofof}[1]{\par\noindent {\it Proof of #1}.}
\newcommand{\Step}[1]{\par\noindent {\it Step} #1.}

\def\Div{\mathop{\mathrm{div}}\nolimits}
\def\rot{\mathop{\mathrm{rot}}\nolimits}
\def\Re{\mathop{\mathrm{Re}}\nolimits}
\def\Im{\mathop{\mathrm{Im}}\nolimits}

\begin{document}


 \begin{center}
{\LARGE Global Existence and Uniqueness of Solutions to \\
the Maxwell-Schr\"odinger Equations}

\bigskip

\bigskip

{\large Makoto NAKAMURA\footnote{%
Supported by Grant-in-Aid for Young Scientists (B) \#16740071 of 
Japan Ministry of Education, Culture, Sports, Science and Technology.}
}


{\small \it
Mathematical Institute\textup{,} Tohoku University \\
Sendai 980-8578\textup{,} Japan \\ \rm
E-mail: makoto@math.tohoku.ac.jp}

\smallskip

and 

\smallskip

{\large Takeshi WADA\footnote{%
Supported by Grant-in-Aid for Young Scientists (B) \#16740075 of 
Japan Ministry of Education, Culture, Sports, Science and Technology.}
}

{\small \it
Department of Mathematics\textup{,} Faculty of Engineering\textup{,} Kumamoto University \\
Kumamoto 860-8555\textup{,} Japan  \\ \rm
E-mail: wada@gpo.kumamoto-u.ac.jp}

\bigskip

\revised{\it Dedicated to Professor Hiroki Tanabe on his seventy-fifth birthday}
\end{center}

\bigskip

\noindent{\bf Abstract.} 
The time local and global well-posedness for the Maxwell-Schr{\"o}dinger equations is considered in Sobolev spaces in three spatial dimensions.
The Strichartz estimates of Koch and Tzvetkov type are used for obtaining the solutions in the Sobolev spaces of low regularities.
One of \revised{the} main results is that the solutions exist time globally for large data.
\bigskip

\section{Introduction}\label{sec:intro}
The Maxwell-Schr\"odinger system  (MS) in space dimension 3 describes the time evolution of 
a charged nonrelativistic quantum mechanical particle interacting with 
the (classical) electro-magnetic field it generates.
We can \revised{state} this system in usual vector notation  as follows:
\begin{gather}
i\partial_tu=(-\DeltaA+\phi)u,   \label{eq:MS1} \\
-\Delta \phi -\partial_t \Div \A =\rho,  \label{eq:MS2} \\
\square \A + \nabla (\partial_t \phi +\Div \A)  =\J, \label{eq:MS3}
\end{gather}
where $(u,\phi,\A): \R^{1+3}\rightarrow \C\times\R\times\R^3$, 
$\nablaA=\nabla -i\A$, 
$\DeltaA={\nablaA}^2$, $\rho=|u|^2$, $\J=2\Im \bar{u}\nablaA u$,
and $\nabla$, $\Delta$ and $\square$ are the usual gradient, Laplacian and d'Alembertian respectively.
Physically, $u$ is the wave function of the particle, $(\phi, \A)$ is the electro-magnetic potential,
$\rho$ is the charge density, and $\J$ is the current density.

The system (MS)  formally conserves at least two quantities, namely the total charge 
$\mathcal{Q}\equiv\| u \|_2^2$
and the total energy
\[\mathcal{E}\equiv \| \nablaA u\|_2^2 +\frac 12 \| \nabla \phi +\partial_t \A \|_2^2 
+ \frac 12\| \rot \A\|_2^2. \]
The system (MS) is invariant under the gauge transform
\begin{equation}
(u',\phi',\A')=(\exp(i\lambda)u, \phi-\partial_t\lambda,\A+\nabla\lambda)
\label{eq:gauge}
\end{equation}
and in this paper we mainly study it in the Coulomb gauge 
\begin{equation}\label{eq:coulomb}
\Div \A=0, 
\end{equation}
in which we can treat the system most easily.
In this gauge, \eqref{eq:MS2} and \eqref{eq:MS3} become
\begin{equation}
-\Delta \phi=\rho,\quad \square \A+ \nabla \partial_t \phi =\J. \label{eq:MS23}
\end{equation}
The first equation of \eqref{eq:MS23} is solved as
\[ \phi=\phi(u)=(-\Delta)^{-1} \rho =(4\pi |x|)^{-1} *|u|^2 \]
and the term $\nabla \partial_t \phi$ in the second equation is dropped by operating 
the Helmholtz projection $P=1-\nabla\Div\Delta^{-1}$ to the both sides of the equation.
Therefore in the Coulomb gauge the system (MS) is rewritten as
\begin{gather}
i\partial_tu=(-\DeltaA+\phi(u))u, \label{eq:MSC1}\\
\square \A =P\J, \label{eq:MSC2}
\end{gather}
which is referred to as (MS-C).
To solve (MS-C) we should give the initial condition
\begin{equation}
(u(0), \A(0), \partial_t \A(0)) =(u_0,\A_0,\A_1) \label{eq:IDC}
\end{equation}
in the direct sum of Sobolev spaces
\[ X^{s,\sigma}= \{(u_0,\A_0,\A_1)\in H^s\oplus H^\sigma\oplus H^{\sigma-1} ;
\Div \A_0=\Div \A_1=0\}.\]
The  condition \eqref{eq:coulomb} is conserved under the consistency conditions $\Div \A_0=\Div \A_1=0$
since the equation  $\square \Div \A=0$ follows from \eqref{eq:MSC2}.

Several authors have studied the Cauchy problem and the scattering theory for (MS-C).
Nakamitsu-M. Tsutsumi \cite{NT86} showed the time local well-posedness for (MS-C) in 
$X^{s,\sigma}$ with $s=\sigma=3,4,5,\dots$.
In fact, they treated the case of Lorentz gauge mentioned below, 
but the Coulomb gauge case can be treated analogously.
We remark that their condition can be  refined  as $s=\sigma>5/2$
by the use of fractional order Sobolev spaces and the commutator estimate by Kato-Ponce~\cite{KP88}.
Recently Nakamura-Wada~\cite{NW05} showed the time local well-posedness for wider class of 
$(s,\sigma)$ including the case $s=\sigma\ge5/3$ (precisely see the remark for 
Theorem~\ref{thm:MS-C}) by using covariant derivative estimates for the Schr\"odinger part 
and the Strichartz estimate for the Maxwell part.
On the other hand, Guo-Nakamitsu-Strauss \cite{GNS95} constructed 
a time global (weak) solution in $X^{1,1}$ although they did not show the uniqueness.
Indeed, in the Coulomb gauge the energy takes the form 
\[\mathcal{E}= \| \nablaA u\|_2^2 +\frac 12 \| \nabla \phi\|_2^2 
+\frac 12 \| \partial_t \A \|_2^2 + \frac 12\| \nabla \A\|_2^2, \]
and hence $\| (u,\A,\partial_t \A);X^{1,1}\|$ does not blow up. 
Therefore the global existence is proved by parabolic regularization and compactness method.
For the scattering theory, the existence of modified wave operators was proved by 
Y. Tsutsumi~\cite{T93}, Shimomura~\cite{S03}, and Ginibre-Velo~\cite{GV03,GVpre}.

As we have summarized above, there are several results for the Cauchy problem 
both at $t=0$ or $t=\infty$.
However there are no results concerning the global existence of strong solutions 
even for small data;
the solutions to
 (MS-C) obtained in~\cite{GV03,GVpre,S03,T93}  exist only for $t\ge 0$ 
 and we do not know whether these solutions globally exist or blow-up at finite negative time.
 The aim of this paper is to answer this problem.
Shortly, we prove the global existence of unique strong solutions.
To do this, we  would need a priori estimates derived from the conservation laws of charge and energy, 
and hence it is desirable to show the local well-posedness in lower regularity.
Therefore we first refine the local theory.
To make the statements of the  propositions simple, we introduce the notation
\begin{align*}
\mathcal{R}_* &=\{ (s,\sigma)\in \R^2; \sigma \ge \max\{1;s-2;(2s-1)/4\}, 
(s,\sigma)\neq (7/2,3/2) \}, \\
\mathcal{R}^* &=\{ (s,\sigma)\in \R^2; \sigma \le \min\{s+1;3s/2;2s-3/4 \}, 
(s,\sigma)\neq (2,3) \} 
\end{align*}
and $\mathcal{R}=\mathcal{R}_* \cap \mathcal{R}^*$.
\begin{theorem}
\label{thm:MS-C}\samepage
Let $(s, \sigma)\in \mathcal{R}$ with
$s\ge 11/8$\textrm{,}  $\sigma>1${\rm.}
Then for any $(u_0,\A_0,\A_1)\in X^{s,\sigma}${\rm ,}
there exists $T>0$ such that \textup{(MS-C)} with initial condition \eqref{eq:IDC} 
has a unique solution $(u,\A)$ satisfying
$(u,\A, \partial_t\A)\in C([0,T]; X^{s,\sigma})${\rm .}
Moreover if $s>11/8$ and $(s+1,\sigma)\in \mathcal{R}_*${\rm ,}
then the mapping $(u_0,\A_0,\A_1) \mapsto (u,\A,\partial_t\A)$ is continuous 
as a mapping from $X^{s,\sigma}$ to $C([0,T]; X^{s,\sigma})${\rm .}
\end{theorem}


\begin{figure}[tbh]
\begin{center}
\begin{picture}(350,350)

\put(20,40){\vector(1,0){280}}
\put(80,40){\circle*{2}}
\put(120,40){\circle*{2}}
\put(160,40){\circle*{2}}
\put(200,40){\circle*{2}}
\put(240,40){\circle*{2}}
\put(280,40){\circle*{2}}

\put(20,22){0}
\put(78,22){1}
\put(118,22){2}
\put(158,22){3}
\put(198,22){4}
\put(238,22){5}
\put(278,22){6}
\put(310,35){$s$}

\put(40,20){\vector(0,1){280}}
\put(40,80){\circle*{2}}
\put(40,120){\circle*{2}}
\put(40,160){\circle*{2}}
\put(40,200){\circle*{2}}
\put(40,240){\circle*{2}}
\put(40,280){\circle*{2}}

\put(26,77){1}
\put(26,117){2}
\put(26,157){3}
\put(26,197){4}
\put(26,237){5}
\put(26,277){6}
\put(35,310){$\sigma$}

\put(95,20){\line(0,1){270}}
\put(80,305){$s=\frac{11}{8}$}
\put(50,30){\line(1,2){130}}
\put(130,305){$\sigma=2s-\frac{3}{4}$}
\put(32,28){\line(2,3){175}}
\put(205,305){$\sigma=\frac{3s}{2}$}
\put(20,60){\line(1,1){230}}
\put(260,305){$\sigma=s+1$}
\put(120,160){\circle{3}}
\put(125,155){$(2,3)$}
\put(100,130){\circle*{3}}
\put(105,125){$(\frac{3}{2},\frac{9}{4})$}
\put(100,20){\line(1,1){190}}
\put(300, 210){$\sigma=s-2$}
\put(30,25){\line(2,1){260}}
\put(300,150){$\sigma=\frac{2s-1}{4}$}

\multiput(30,45)(4,2){65}{.}
\put(300,180){$\sigma=\frac{2s+1}{4}$}
\multiput(60,20)(3,3){75}{.}
\put(300, 240){$\sigma=s-1$}
\put(142,101){\circle{3}}
\put(108,105){$(\frac{5}{2},\frac{3}{2})$}
\put(180,100){\circle{3}}
\put(185,90){$(\frac{7}{2},\frac{3}{2})$}

\multiput(33,80)(4,0){64}{.}
\put(300,77){$\sigma=1$}

\put(90,-15){A figure of the range of $s$ and $\sigma$.}

\end{picture}
\end{center}
\bigskip
\end{figure}


\begin{remark}
(1)
$T$ depends only on $s,\sigma$ and $\|(u_0,\A_0,\A_1);X^{s,\sigma}\|$.

(2)
For any $s$ and $\sigma$ satisfying the assumption above for the unique existence of the solution, 
the mapping $(u_0,\A_0,\A_1) \mapsto (u,\A,\partial_t\A)$ 
is continuous in w*-sense.
Namely if a sequence of initial data strongly converges in $X^{s,\sigma}$, 
then corresponding sequence of solutions also converges star-weakly in $L^\infty(0,T; X^{s,\sigma})$.

(3) In \cite{NW05}, we also assume $s\ge 5/3$ and $4/3\le\sigma \le (5s-2)/3$
with $(s,\sigma)\neq (5/2,7/2)$.
\end{remark}

\revised{%
Generally, in order to construct solutions of dispersive equations in low regularity function spaces,
we usually use smoothing effects such as Strichartz estimates.
However, usual Strichartz estimates for Schr\"odinger equations does not match the equation \eqref{eq:MSC1}
since we cannot avoid the loss of derivative coming from the term $2i\A\cdot \nabla u$.
This is why the preceding results rely on the $L^2$-based energy method.
In the present work we use a variation of Strichartz estimates first given by 
Koch-Tzvetkov~\cite{KT03} and refined by Kenig-Koenig~\cite{KK03}
for Benjamin-Ono type equations, and adapted for Schr\"odinger equations by J. Kato~\cite{K05}.
In the proof of Theorem~\ref{thm:MS-C}, we slightly refine this estimate and combine it 
with the covariant derivative estimates developed in our previous work~\cite{NW05}. 
Our local theory does not cover the result for the energy class $H^1$,
but it is sufficient for our aim.
Indeed, we can show the following global result:
}
\begin{theorem}\label{thm:global}
The solution obtained in Theorem {\rm\ref{thm:MS-C}} exists time globally{\rm .}
\end{theorem}

\revised{By the use of Koch-Tzvetkov type estimate, we can show 
\[ \| u; L^2(0,T; H^{1/2-\delta}_6)\|\le C \langle T \rangle^3, \]
where $\delta>0$ is sufficiently small and the constant $C$ depends only on 
$\| (u_0,\A_0,\A_1);X^{1,1} \|$ and $\delta$.
Roughly speaking, we can gain $1/2$ regularity by this estimate.
Indeed, if we control this norm by the Sobolev inequality, we would need $\| u;H^{3/2-\delta}\|$.
This estimate, combined with several a priori estimates obtained in the proof of Theorem~\ref{thm:MS-C},
shows that $\|(u,\A,\partial_t \A);X^{s,\sigma}\|$ does not blow up.
}

Next we consider the Lorentz gauge 
\begin{equation}
\partial_t\phi+\Div \A=0.
\label{eq:lorentz}
\end{equation}
(MS)  in the Lorentz gauge, which is referred to as (MS-L),  is expressed as
\begin{gather*}
i\partial_tu=(-\DeltaA+\phi)u, \quad 
\square \phi =\rho, \quad 
\square\A =\J.
\end{gather*}
In this case, we need the initial data
\begin{equation}\label{eq:IDL}
(u(0),\phi(0),\partial_t \phi(0), \A(0), \partial_t \A(0))
=(u_0,\phi_0,\phi_1,\A_0,\A_1)\in Y^{s,\sigma}.
\end{equation}
Here
\begin{align*}
Y^{s,\sigma}&=\{ (u_0,\phi_0,\phi_1,\A_0,\A_1) \in 
H^s\oplus H^\sigma\oplus H^{\sigma-1}\oplus H^\sigma\oplus H^{\sigma-1}; \\
&\quad \Div \A_0+\phi_1=\Div \A_1+\Delta \phi_0+|u_0|^2=0 \}. 
\end{align*}
The condition \eqref{eq:lorentz} is conserved under the consistency condition in the definition of  $Y^{s,\sigma}$
since $\square (\partial_t\phi+\Div \A)=\partial_t\rho +\Div \J=0$.
The first and the  second equations respectively follow from the wave equations both for $\phi$ and $\A$, 
and from the conservation of charge derived from the Schr\"odinger equation.
Our result for (MS-L)  is the following.

\begin{theorem}\samepage
\label{thm:MS-L}
Let $(s,\sigma)\in \mathcal{R}^*$ with $s\ge 11/8, \sigma>1$ and $\sigma \ge s-1${\rm .}
Then for any
$(u_0,\phi_0,\phi_1,\A_0,\A_1)\in Y^{s,\sigma}${\rm ,}
there exists $T>0$ such that 
\textup{(MS-L)} with initial condition \eqref{eq:IDL} has a unique solution 
$(u,\phi,\A)$ satisfying 
\[(u,\phi,\partial_t \phi,\A,\partial_t \A) \in C([0,T];Y^{s,\sigma}). \] 
This solution exists time globally\textup{.}
Moreover{\rm ,} if $s>11/8$ and  
$(s+1,\sigma)\in \mathcal{R}_*${\rm ,} 
then the mapping
$(u_0,\phi_0,\phi_1,\A_0,\A_1)\mapsto (u,\phi,\partial_t \phi,\A,\partial_t \A)$ 
is continuous as a mapping from $Y^{s,\sigma}$ to $C([0,T]; Y^{s,\sigma})${\rm .}
\end{theorem}


We can also consider the temporal gauge 
\begin{equation}
\phi=0.
\label{eq:temporal}
\end{equation}
(MS)  in the temporal gauge, which is referred to as (MS-T),  is expressed as
\begin{gather*}
i\partial_tu=-\DeltaA u, \quad 
\square\A +\nabla {\mbox {\rm div}} \A=\J(u,\A).
\end{gather*}
In this case, we need the initial data
\begin{equation}\label{eq:IDT}
(u(0), \A(0), \partial_t \A(0))
=(u_0, \A_0,\A_1)\in Z^{s,\sigma}.
\end{equation}
Here
\[
Z^{s,\sigma}=\{ (u_0,\A_0,\A_1) \in 
H^s\oplus H^\sigma\oplus H^{\sigma-1} \  ; \ -\Div \A_1=|u_0|^2 \}. 
\]
Our result for (MS-T)  is the following.

\begin{theorem}\samepage
\label{thm:MS-T}
Let $(s,\sigma)\in \mathcal{R}^*$ with $s\ge 11/8, \sigma>1$ and $\sigma \ge s-1${\rm .}
Then for any
$(u_0,\A_0,\A_1)\in Z^{s,\sigma}${\rm ,}
there exists $T>0$ such that 
\textup{(MS-T)} with initial condition \eqref{eq:IDT} has a unique solution 
$(u,\A)$ satisfying 
\[(u,\A,\partial_t \A) \in C([0,T];Z^{s,\sigma}). \] 
This solution exists time globally\textup{.}
Moreover{\rm ,} if $s>11/8$ and  
$(s+1,\sigma)\in \mathcal{R}_*${\rm ,} 
then the mapping
$(u_0,\A_0,\A_1)\mapsto (u,\A,\partial_t \A)$ 
is continuous as a mapping from $Z^{s,\sigma}$ to $C([0,T]; Z^{s,\sigma})${\rm .}
\end{theorem}


This paper is organized as follows:
In \S 2,  We first prepare basic estimates used throughout this paper,
namely an estimate of Hartree type nonlinearities (Lemma~\ref{lem:phiest}) and 
that for covariant derivatives (Lemma~\ref{lem:cov}).
Next we introduce Strichartz estimates for Klein-Gordon equations (Lemma~\ref{lem:StKG})
and Koch-Tzvetkov type Strichartz estimates for Schr\"odinger equations 
(Lemma \ref{lem:KT}).
We also prepare the estimate of the nonlinear term of the Maxwell part\revised{, which is} based on the Kato-Ponce commutator estimate 
(Lemmas \ref{lem:commutator} and \ref{lem:nlmaxwell}).
In \S 3, we study the linearized Schr\"odinger equation associated with \eqref{eq:MSC1}.
Applying the Koch-Tzvetkov type estimate, we derive a smoothing property of 
the Schr\"odinger equation with electro-magnetic potential  (Lemma~\ref{lem:smoothing}). 
Using this estimate together with covariant derivative estimates, 
we prove the unique solvability of this equation first in $H^2$ 
(Lemma~\ref{lem:H2est})
and next in $H^s$ with $s\ge 0$ (Lemmas~\ref{lem:propertyU} and~\ref{lem:Hsest}).
In \S 4, we discuss the local solvability.
We prove the local well-posedness by the contraction mapping principle 
(Propositions~\ref{prop:localexistence} and~\ref{prop:regularity}).
The continuous dependence of the solutions on the data is left to \S 6, 
since it is usually the most delicate part of the theory of well-posedness.
In \S 5, we derive a priori estimates of solutions (Lemma~\ref{lem:apriori}) 
and use them in the proof of global existence.
\S 6 is devoted to the proof of the continuous dependence of the solutions on the data.
In \S 7, we prove Theorems~\ref{thm:MS-L} and \ref{thm:MS-T} by using the gauge transform.

We conclude this section by giving the notation used in this paper.
$L^p=L^p(\R^3)$ is the usual Lebesgue space 
and its norm is denoted by $\|\cdot\|_p$.
$p'=p/(p-1)$ is the dual exponent of $p$. 
This symbol is used only for Lebesgue exponents.
$H^s_p=\{ \phi \in \mathcal{S}'(\R^3); \| (1-\Delta)^{s/2} \phi \|_p <\infty \}$ 
is the usual Sobolev space.
$\Dot{H}^s_p=\{ \phi \in \mathcal{S}'(\R^3); \| (-\Delta)^{s/2} \phi \|_p <\infty \}$ 
is the homogeneous Sobolev space.
\revised{The subscript $p$ is omitted if $p=2$.}
For any interval $I \subset \R$ and Banach space $X$, $L^p(I;X)$ denotes the space of 
$X$-valued strongly measurable functions on $I$ whose $X$-norm belong to $L^p(I)$.
This space is often abbreviated to $L_T^p X$ for  $I=(0,T)$.
Similarly we use the abbreviation $C_T^m X=C^m ([0,T];X)$ and 
$W_T^{m,p}X=W^{m,p}(0,T;X)$,
where $W^{m,p}(I;X)$ denotes the space of functions in $L^p(I;X)$ 
whose derivatives up to the $(m-1)$-times are locally absolutely continuous 
and the derivatives up to the $m$-times belong to $L^p(I;X)$.
For normed spaces $X_i, i=1,\dots,n$, we define the norm of $X=\bigcap_{i=1}^n X_i$ by
$\| \cdot; X \| =\max_{i=1}^n \| \cdot;X_i \|$ so that $X$ is also a normed space.
We define 
$\Sigma^{m,s}_T=\bigcap_{j=0}^m W^{j,\infty}_T H^{s-2j}$
and
$\M^{m,\sigma}_T=\bigcap_{j=0}^m W^{j,\infty}_T H^{\sigma-j}$.
The inequality $a\lesssim b$ means $a\le Cb$, 
where $C$ is a positive constant that is not essential.
We write $a\simeq b$ if $b\lesssim a$ as well as $a\lesssim b$.
$\langle a\rangle=\sqrt{1+a^2}$. 
$a\vee b$ and $a\wedge b$ denote the maximum and the minimum of $a$ and $b$ respectively.
We use the following unusual but convenient symbol:
$a_+$ means $a\vee0$ if $a\neq0$, whereas $0_+$ means a sufficiently small positive number.
Namely
$b\ge a_+$ means $b\ge a\vee 0$ if $a\neq 0$, and $b>0$ if $a=0$.
It is useful to express 
sufficient conditions for Sobolev type embeddings $H^s_r \embedding L^p$
by the inequality $(1/r-s/3)_+ \le 1/p\le 1/r$ with $1\le r<\infty$.


\section{Preliminaries}\label{sec:pre}
In this section we summarize lemmas used in the proof of 
Theorems~\ref{thm:MS-C}-\ref{thm:MS-T}.
The following two lemmas will be repeatedly used in estimates of nonlinear terms.
\begin{lemma} \samepage
\label{lem:phiest}
Let $s, s_1, s_2, s_3$ be nonnegative numbers satisfying $s\le s_3$ and 
$s_1\wedge s_2\ge s-2$ with $s_1+s_2>0${\rm .}
Let 
$s_1+s_2+s_3\wedge (3/2)\ge s+1$
and the inequality be strict if {\rm (i)} $s_j=3/2$ for some $1\le j\le 3$
or {\rm (ii)} $s=s_3<3/2${\rm .}
Then the following estimate holds{\rm:}
\begin{equation}\label{eq:phiest}
\| (-\Delta)^{-1} (u_1 u_2) u_3 ;H^s \| 
\lesssim 
\prod_{j=1}^3 \| u_j ; H^{s_j} \|. 
\end{equation}
\end{lemma}
\Proof
See Lemma 2.1 in ~\cite{NW05}.

\begin{lemma} \samepage \label{lem:cov}
Let $(s,\sigma)\in \mathcal{R}_*$ with $s\ge 0$\textup{.} 
Let $\A\in H^\sigma$ satisfy $\Div \A=0${\rm .}
\par \indent
\textup{(i)}
Let $V_1(\A,v)=2i \A\cdot \nabla v +|\A|^2 v$\textup{.}
Then $V_1$ is a continuous mapping from $H^\sigma \times H^{s-1/2}$ to $H^{s-2}$
with the estimate
\begin{equation} \label{eq:V1}
\| V_1(\A,v);H^{s-2}\| \lesssim \langle \| A;H^\sigma \| \rangle^2 \| v;H^{s-1/2}\|.
\end{equation}
\par \indent
\textup{(ii)} 
The following estimates hold for any $v\in H^s$\textup{:}
\begin{align}
\|v;H^s\| +\langle\|\A;H^\sigma\| \rangle^\alpha \|v\|_2 \simeq \| \DeltaA v;H^{s-2}\|
+\langle\|\A;H^\sigma\| \rangle^\alpha \|v\|_2,
\label{eq:cov}
\end{align}
where $\alpha=\alpha(s,\sigma)$ is a positive constant independent of $v$ and $\A${\rm .}
Especially if $s=2$\textup{,} the estimate
\begin{equation}\label{eq:cov2}
\|v;H^2\| +\langle\|\A;\Dot{H}^1\| \rangle^4 \|v\|_2  
\simeq \| \DeltaA v\|_2+\langle\|\A;\Dot{H}^1\| \rangle^4 \|v\|_2
\end{equation}
holds valid\textup{.}
\end{lemma}

\Proof
For $s\ge 5/2$, we can prove \eqref{eq:V1} by the Leibniz rule and the Sobolev inequality.
The case $s=0$ is the dual of the case $s=5/2$, and hence the case $0<s<5/2$ is proved by interpolation.
The continuity of $V_1$ is proved in the same way.
Applying the interpolation inequality $\| v;H^{s-1/2}\|\le \| v;H^s\|^{1-1/2s} \|v\|_2^{1/2s}$ 
to the estimate of $V_1$, we can show \eqref{eq:cov}. We can show \eqref{eq:cov2}
if we use the estimate $\| V_1 \|_2 \le \| \A \|_6 \| \nabla v\|_3+\| \A\|_6^2 \|v\|_6$ 
instead of \eqref{eq:V1}. \QED

\thmskip
Next we introduce Strichartz type estimates  for Klein-Gordon equations
(see for example~\cite{B84,GV85AIHP,GV95,S77}).
\begin{lemma}\label{lem:StKG} \samepage
Let $T>0$, $\sigma\in \R$ and let $(q_j,r_j)$\textup{,} $j=0,1$\textup{,} 
satisfy $0\le 2/q_j=1-2/r_j<1$\textup{.}
Let $(A_0,A_1)\in H^\sigma \oplus H^{\sigma-1}$ and 
$F\in L_T^{{q_1}'}H^{\sigma-1+2/q_1}_{{r_1}'}$\textup{.}
Then a solution $A$ to the equation 
\[ (\square+1) A=F \]
 with $A(0)=A_0, \partial_t A(0) =A_1$ belongs to
$C_T H^\sigma \cap C^1_T H^{\sigma-1}$ and satisfies the estimate
\begin{equation}
\max_{k=0,1}\| \partial_t^k A;L_T^{q_0}H^{\sigma-k-2/q_0}_{r_0}\|
\lesssim 
\| (A_0, A_1);H^\sigma \oplus H^{\sigma-1}\|+\| F;L_T^{{q_1}'}H^{\sigma-1+2/q_1}_{{r_1}'} \|.
\end{equation}
\end{lemma}

\begin{lemma} \label{lem:KT} \samepage
Let $T>0$\textup{,} $s\in \R$\textup{,} $\alpha>0$ and $0\le 2/q=3/2-3/r\le 1$\textup{.}
Let $f \in L_T^2H^{s-2\alpha}$.
Then a solution $u\in \Sigma^{1,s}_T$ to the equation
\[i\partial_t u=-\Delta u +f \]
belongs to $L_T^q H^{s-\alpha}_r$ and satisfies the estimate
\begin{equation} \label{eq:KT}
\| u;L_T^q H^{s-\alpha}_r\|
\lesssim \| u ; L_T^\infty H^s\|
+T^{1/2}\| f;L_T^2H^{s-2\alpha} \|.
\end{equation}
\end{lemma}

\begin{remark}
This kind of estimates was first given by Koch-Tzvetkov~\cite{KT03} 
\revised{for the Benjamin-Ono equation} , and it is  Kenig-Koenig~\cite{KK03}
who formulated the estimate as above. 
\revised{Kato~\cite{K05} adapted this estimate for Schr\"odinger equations.}
However, in \cite{KK03,K05}, they  need an extra assumption $u\in L_T^\infty H^{s+\epsilon}$ 
to prove \eqref{eq:KT}, 
with the first term in the right-hand side replaced by $\| u; L_T^\infty H^{s+\epsilon}\|$.
\end{remark}
\Proof
Without loss of generality we assume $s=0$.
Let $u=\sum_{j=0}^\infty u_j$ is the Littlewood-Paley decomposition of the solution $u$.
Namely we take $\psi\in\mathcal{S}'(\R^3)$ such that 
$\textrm{supp}\Hat{\psi}\subset\{\xi; 1/2\le | \xi | \le 2 \}$ and
$\sum_{j=-\infty}^\infty \Hat{\psi}(\xi/2^j)=1$ for any $\xi \neq 0$,
and put $u_j=\mathcal{F}^{-1} \left(\Hat{\psi}(\xi/2^j)\Hat{u}(t,\xi)\right)$ for $j\ge 1$,
$u_0=\mathcal{F}^{-1} \sum_{j=-\infty}^0\Hat{\psi}(\xi/2^j)\Hat{u}(t,\xi)$.
Here $\Hat{u}(t,\xi)$ is the Fourier transform with respect to the space variable.
Similarly let $f=\sum_{j=0}^\infty f_j$ is the Littlewood-Paley decomposition of  $f$.
Then $u_j$ satisfy the equation
\begin{equation}\label{eq:decomposed}
i\partial_t u_j=-\Delta u_j +f_j.
\end{equation}
We divide the interval $[0,T]$ into disjoint intervals $\{I_{jk}\}_{k=1}^{m_j}$
such that $2^{-2\alpha j}T\le |I_{jk}|<2^{-2\alpha j+1}T$ and take $t_{jk}\in \Bar{I}_{jk}$
at which $\| u_j(t)\|_2$ attains its minimum in the interval $\Bar{I}_{jk}$.
By the standard Strichartz estimate for Schr\"odinger equations~\cite{K87,KT98,S77,Y87},
\[  \| u_j;L^q(I_{jk}; L^r)\| \lesssim \| u_j(t_{jk})\|_2+ \| f_j;L^1(I_{jk}; L^2)\|. \]
Taking the sum with respect to $k$, we obtain
\begin{align*}
\| u_j;L^q_T L^r\| &=(\sum_{k=1}^{m_j} \| u_j;L^q(I_{jk}; L^r)\|^q)^{1/q}
\le (\sum_{k=1}^{m_j} \| u_j;L^q(I_{jk}; L^r)\|^2)^{1/2} \\
&\lesssim (\sum_{k=1}^{m_j} \| u_j(t_{jk})\|_2^2)^{1/2}
+(\sum_{k=1}^{m_j} \| f_j;L^1(I_{jk}; L^2)\|^2)^{1/2} \\
&\le (\sum_{k=1}^{m_j} T^{-1} 2^{2\alpha j} |I_{jk}| \| u_j(t_{jk})\|_2^2)^{1/2}
+(\sum_{k=1}^{m_j} 2^{-2\alpha j+1}T\| f_j;L^2(I_{jk}; L^2)\|^2)^{1/2} \\
&\le T^{-1/2} \| 2^{\alpha j}u_j; L_T^2L^2 \|
+(2T)^{1/2}\|  2^{-\alpha j}f_j;L_T^2L^2 \|.
\end{align*}
The first term in the right-hand side is obtained from the fact that $t_{jk}$ are the minimum points
and the definition of integral.
We have also used the definition of $I_{jk}$and the H\"older inequality for the time variable.
Therefore
\begin{align*}
\| u ;L^q_T H^{-\alpha}_r\| 
&\simeq \| \{\sum_{j=0}^\infty (2^{-\alpha j}|u_j |)^2\}^{1/2};L^q_T L^r\| 
\le  (\sum_{j=0}^\infty\|  2^{-\alpha j}|u_j | ;L^q_T L^r\|^2)^{1/2} \\
&\lesssim T^{-1/2} (\sum_{j=0}^\infty\| u_j; L_T^2L^2 \|^2)^{1/2} 
+T^{1/2}(\sum_{j=0}^\infty\|  2^{-2\alpha j}f_j;L_T^2L^2 \|^2)^{1/2} \\
&\simeq T^{-1/2}\| u; L_T^2L^2 \|
+T^{1/2}\|  f ;L_T^2H^{-2\alpha} \|,
\end{align*}
where we have used the equivalent norms between the Sobolev spaces and the Triebel-Lizorkin spaces (see page 29 in \cite{Triebel92}).
Thus the lemma has been proved.
\QED

 \begin{lemma}\label{lem:commutator} \samepage
Let $\sigma\ge 0$\textup{.} 
Let $1<p, p_1, p_4<\infty$ and $1<p_2, p_3\le \infty$ satisfy $1/p=1/p_1+1/p_2=1/p_3+1/p_4$\textup{.}
Then the following estimate holds valid\textup{:}
\begin{equation}\label{eq:commutator}
\| P(\Bar{u}_1\nabla u_2); H^\sigma_p\|
\lesssim \|u_1 ; H^\sigma_{p_1}\| \| \nabla u_2 \|_{p_2}
+\| \nabla u_1 \|_{p_3} \|u_2 ; H^\sigma_{p_4}\|.
\end{equation}
Moreover if $0\le \sigma \le 1$\textup{,} we can omit the second term of the right-hand side\textup{.}
\end{lemma}
\Proof
We can prove \eqref{eq:commutator} by the use of the Kato-Ponce commutator estimate~\cite{KP88}
and the fact that $P\nabla=0$; 
see~\cite{NW05} for detail.
The last assertion can be checked immediately for $\sigma=0,1$ and generalized for $0<\sigma<1$
by interpolation.
\QED

\begin{lemma}\label{lem:nlmaxwell}\samepage
Let $(s,\sigma)\in \mathcal{R}^*$ with $s\ge 5/4$ and  $\sigma\ge 1$\textup{.}
Then 
\begin{equation}
\| P(\Bar{u}_1\nablaA u_2 ); L^1_T H^{\sigma-1}\|\lesssim
T^{1/4}\langle T \rangle^{3/4} \langle \| \A; L^\infty_T H^{\sigma-1/2}\| \rangle
\prod_{j=1}^2\| u_j ;L^\infty_T H^s \cap L^2_T H^{s-1/2}_6 \|.
\end{equation}
\end{lemma}
\Proof
It suffices to show the following inequalities:
\begin{align}
\| P(\Bar{u}_1\nabla u_2); L^1_T H^{\sigma-1}\| &\lesssim
T^{1-1/q_1-1/q_2}
(\| u_1 ;L^{q_1}_T H^{s-1/q_1}_{r_1} \|\| u_2 ;L^{q_2}_T H^{s-1/q_2}_{r_2} \| \nonumber \\
&\qquad+\| u_1 ;L^{q_2}_T H^{s-1/q_2}_{r_2}\|\| u_2 ;L^{q_1}_T H^{s-1/q_1}_{r_1} \|) 
\label{eq:subestimate1}
\end{align}
for suitable $q_j, r_j$ satisfying the conditions $0\le 2/q_j=3/2-3/r_j\le 1$ and
$2/q_1+2/q_2\le 3/2$;
\begin{equation}\label{eq:subestimate2}
\| \A\Bar{u}_1 u_2; L^1_T H^{\sigma-1}\|\lesssim
T\|\A; L^\infty_T H^{\sigma-1/2}\|  \prod_{j=1}^2\| u_j ;L^\infty_T H^s \|.
\end{equation}
We first prove \eqref{eq:subestimate1}. We use Lemma~\ref{lem:commutator} and obtain the estimate
\[ \| P(\Bar{u}_1\nabla u_2); H^{\sigma-1}\|
\lesssim \|u_1 ; H^{\sigma-1}_{p_1}\| \| \nabla u_2 \|_{p_2}
+\| \nabla u_1 \|_{p_2} \|u_2 ; H^{\sigma-1}_{p_1}\|, \]
where the choice of $p_1, p_2$ depends on the value of $s$.
Practically, if $s>2$ and $\sigma \le s+1$, we choose $p_1=2, p_2=\infty$ so that 
$H^s\embedding H^{\sigma-1}$ and $H^{s-3/2}_6 \embedding L^\infty$.
Then we obtain \eqref{eq:subestimate1} with $2/q_1=0, 2/q_2=1$.
If $3/2\le s<2$, we choose $1/p_1=s/3-1/6, 1/p_2=(2-s)/3$
so that $H^{s-1/q_2}_{r_2}\equiv H^{s-1/2}_6\embedding H^1_{p_2}$ by the Sobolev inequality.
Putting $r_1=p_1$ and hence $2/q_1=2-s$, we obtain 
$H^{s-1/q_1}_{r_1}\embedding H^{\sigma-1}_{p_1}$ and $2/q_1+2/q_2=3-s\le 3/2$
under the condition $\sigma \le 3s/2$.
If $5/4\le s<3/2$, we choose $1/p_1=2(s-1)/3, 1/p_2=1/2-2(s-1)/3$
so that $H^{s-1/q_2}_{r_2}\equiv H^1_{p_2}$.
Putting $r_1=p_1$ and hence $2/q_1=7/2-2s$, we obtain 
$H^{s-1/q_1}_{r_1}\embedding H^{\sigma-1}_{p_1}$ and $2/q_1+2/q_2=3/2$
under the condition $\sigma \le 2s-3/4$.
Therefore we obtain \eqref{eq:subestimate1} by the H\"older inequality for the time variable.
The proof for the case $s=2$ has been omitted, but this  is covered by the case $s \approx 2$.
Next we prove \eqref{eq:subestimate2}.
By the Leibniz rule, 
\begin{align*}
\| \A\Bar{u}_1 u_2; H^{\sigma-1}\| &\lesssim
\|\A; H^{\sigma-1}_3\|  \| u_1\|_{12}\| u_2\|_{12}
+\| \A\|_{p_3}\| u_1;H^{\sigma-1}_{p_4}\| \| u_2\|_{p_5} \\
&\quad+\| \A\|_{p_3}\| u_1\|_{p_5}\| u_2;H^{\sigma-1}_{p_4}\| , 
\end{align*}
where $2\le p_3, p_4, p_5<\infty$ are numbers satisfying
$1/p_3\ge \{(2-\sigma)/3\}_+, 1/p_4 \ge \{ 1/2-(s-\sigma+1)/3\}_+, 1/p_5 \ge (1/2-s/3)_+$
and $1/p_3+1/p_4+1/p_5=1/2$.
Such a choice is possible under the assumption.
Thus \eqref{eq:subestimate2} follows from the Sobolev inequality.
\QED


\section{Linearized Schr\"odinger equations}\label{sec:LS}
In this section we prove some existence theorems and a priori estimates for the linear Schr\"odinger equation associated with \eqref{eq:MSC1}:
\begin{align}\label{eq:LS}
i\partial_t v &= (-\DeltaA +\phi)v, \quad 0<t<T, \\
\label{eq:LSinitial}
v(t_0) &=v_0.
\end{align}
Here $\phi=\phi(u)\equiv (4\pi |x|)^{-1} *|u|^2$ and  $0\le t_0\le T$.
In this section we regard $\A$ and $u$ as known functions defined on $[0,T]$.
Here we clarify the notion of solutions to \eqref{eq:LS}.
\begin{definition}
A function $v$ is called a weak  (resp. strong) $H^s$-solution to \eqref{eq:LS}
if $v$ belongs to $\Sigma^{1,s}_T$ (resp. $C_T H^s \cap C^1_T H^{s-2}$)
and satisfies \eqref{eq:LS} for almost every (resp. all) $0<t<T$.
A function $v$ is called a weak (resp. strong) $H^s$-solution to \eqref{eq:LS}-\eqref{eq:LSinitial}
if $v$ is a weak (resp. strong) $H^s$-solution to \eqref{eq:LS} and satisfies \eqref{eq:LSinitial}.
\end{definition}
For given $\sigma>1$, we fix positive numbers $\delta, q$ and $r$ as
\begin{equation}\label{eq:dqr}
\delta =\min\{(\sigma-1)/2, 1/4\}, \quad 1/q=1/2-2\delta/3, \quad 1/r=2\delta/3.
\end{equation}
This notation will be used throughout the paper.

%
\begin{lemma}\label{lem:smoothing}\samepage
Let $s\ge 0$\textup{,} $\sigma>1$ and  $\sigma\ge s-1$\textup{.}
Let $\A \in L^\infty_T H^\sigma \cap L^2_T L^\infty$ satisfy $\Div \A=0$\textup{.}
Let $u \in L_T^\infty H^{(s-1)\vee 1}$ and $f\in L^2_T H^{s-1}$\textup{.}
Then a weak $H^s$-solution $v$ to 
\[ i\partial_t v = (-\DeltaA +\phi)v+f, \quad 0<t<T, \]
belongs to $L_T^2 H^{s-1/2}_6$ and satisfies the estimate
\begin{align}
\| v; L_T^2 H^{s-1/2}_6\| &\lesssim \langle T \rangle^m  
\langle \| \A;L_T^\infty H^\sigma \cap L^2_T L^\infty \| 
\vee \|u ; L_T^\infty H^{(s-1)\vee 1}\| \rangle^m 
\| v;L_T^\infty H^s\| \nonumber \\
&\quad+T^{1/2} \| f;L^2_T H^{s-1}\|, \label{eq:smoothing}
\end{align}
where $m$ is a positive number\textup{.}
\end{lemma}
\Proof
Applying Lemma \ref{lem:KT},  we obtain
\begin{equation}\label{eq:vL2H6}
\| v; L_T^2 H^{s-1/2}_6\| \lesssim  \| v;L_T^\infty H^s\|
+T^{1/2} \| 2i\A\cdot \nabla v +|\A |^2 v +\phi v+f;L_T^2 H^{s-1}\|. 
\end{equation}
We shall estimate the second term of the right-hand side. We first estimate $2i\A\cdot \nabla v$.
We have 
\begin{equation}\label{eq:leibniz}
\|\A\cdot \nabla v;H^{s-1}\| \lesssim \| \A\|_\infty \| v;H^s\|+
\|\A;H^{s-1}_{r_1}\| \|\nabla v\|_{q_1}
\end{equation}
for $0<2/r_1=1-2/q_1\le 1$, and the second term can be omitted in the case $s\le 1$.
Indeed \eqref{eq:leibniz} without the second term is clearly valid for $s=0$ and $s=1$, 
and it is also valid for $0<s<1$ by interpolation;
on the other hand if $s>1$, \eqref{eq:leibniz} is obtained by the Leibniz rule.
We define the numbers  $q_1, r_1$ as 
(i) $2/q_1=1-2/r_1 =2-s$ if $1<s<2$;
(ii)  $2/q_1=1-2/r_1 =\delta$ if $s=2$;
(iii) $q_1=\infty, r_1 =2$ if $s>2$.
In any case,  by virtue of the Sobolev (Gagliardo-Nirenberg) inequality, 
there exists an exponent $0<\theta <1$ with $2/q_1\le \theta$ such that 
$\| \nabla v\|_{q_1}\lesssim \| v;H^s\|^\theta \| v;H^{s-1/2}_6 \|^{1-\theta}$.
Practically, we can take $\theta=2-s$ if $s<2$, $\theta =3\delta$ if $s=2$,
and $\theta=2(s-2)/(2s-3)$ if $s>2$.
Therefore by the H\"older inequality for the time variable together with the Young inequality,
\begin{align*}
&T^{1/2} \| \A\cdot \nabla v;L_T^2 H^{s-1}\|  
\\
&\quad\lesssim T^{1/2} \|\A; L^2_T L^\infty \| \| v;L_T^\infty H^s\|  
+T^{1/2}\| v; L_T^2 H^{s-1/2}_6\|^{1-\theta}
	\| \A ; L^{2/\theta}_T H^{s-1}_{r_1} \| \| v;L_T^\infty H^s\|^\theta 
\\
&\quad\lesssim \langle T \rangle^m \{ 
 \|\A; L^2_T L^\infty \|  
+\epsilon^{1-1/\theta} \| \A;L_T^{q_1} H^{s-1}_{r_1} \|^{1/\theta}
 \}
	\| v;L_T^\infty H^s\|  
+ \epsilon \| v; L_T^2 H^{s-1/2}_6\|,
\end{align*}
where $m$ is a positive number.
We choose $\epsilon>0$ so small that the last term in the right-hand side is absorbed 
in the left-hand side of \eqref{eq:vL2H6}.
Next we show the estimate
\[ \| |\A |^2 v;H^{s-1}\| \lesssim \| \A;H^\sigma\| \|\A \|_\infty \| v; H^s\|. \]
If $s\le 1$, this inequality follows from the Sobolev inequality. 
If $s>1$, we use the Leibniz rule to derive 
\[ \| |\A |^2 v;H^{s-1}\| \lesssim \| \A\|_3 \|\A \|_\infty \| v; H^{s-1}_6 \|
+ \| \A; H^{s-1}_{r_2} \|  \|\A \|_\infty \| v\|_{q_2} \]
with $0<2/r_2=1-2/q_2\le 1$.
We can choose $q_2$, $r_2$ such that 
$1/r_2\ge (1/2-(\sigma-s+1)/3)_+$ and $1/q_2 \ge (1/2-s/3)_+$.
Therefore we obtain the desired estimate again by the Sobolev inequality.
By the Sobolev inequality $L^{6/(5-2s)} \hookrightarrow H^{s-1}$ for $0\le s\le 1$, and Lemma \ref{lem:phiest} for $s>1$, we have 
\[ \| \phi v;H^{s-1}\|\lesssim \| u;H^{(s-1)\vee 1}\|^2 \| v;H^s\|. \]
Collecting these estimates, we obtain
\begin{align*}
\| v; L_T^2 H^{s-1/2}_6\|  &\lesssim   
\langle T \rangle^m  
\langle \| \A;L_T^\infty H^\sigma \cap L^2_T L^\infty \cap L^{q_1}_T H^{s-1}_{r_1} \| 
\vee \|u ; L_T^\infty H^{(s-1)\vee 1}\| \rangle^m  \\
&\qquad \times \| v;L_T^\infty H^s\| \nonumber \\
&\quad+T^{1/2} \| f;L^2_T H^{s-1}\|, 
\end{align*}
where $m$ is a positive number, and the third space in the norm of $\A$ can be dropped if $s\le 1$.
To complete the proof, we should show 
$\| \A; L^{q_1}_T H^{s-1}_{r_1} \|\lesssim \| \A; L_T^\infty H^\sigma \cap L^2_T L^\infty\|$
in the nontrivial case $1<s\le 2$.
We concentrate on the case $1<s<2$; we can analogously treat the case $s=2$.
Let $\A =\sum_{j=0}^\infty \A_j$ be the Littlewood-Paley decomposition.
Then by the H\"older inequalities both for the sequence and the space variable,
\begin{align*}
\| \A; H^{s-1}_{r_1} \|  
&\simeq \| \{\sum_{j=0}^\infty |2^{(s-1)j} \A_j |^2 \}^{1/2} \|_{r_1} \\
&\lesssim \| \{\sum_{j=0}^\infty |2^{\sigma j }\A_j |^2 \}^{1/2} \|_2^{s-1}
\| \{\sum_{j=0}^\infty |2^{-\epsilon_1 j} \A_j |^2 \}^{1/2} \|_\infty^{2-s},
\end{align*}
where $\epsilon_1=(\sigma-1)(s-1)/(2-s)>0$. 
We have 
$\| \{\sum_{j=0}^\infty |2^{-\epsilon_1 j} \A_j |^2 \}^{1/2} \|_\infty 
\le C_{\epsilon_1}  \| \A\|_\infty$
since $\| \A_j\|_\infty \lesssim \| \A\|_\infty$.
Therefore we obtain
\[ \| \A; L^{q_1}_T H^{s-1}_{r_1} \|\lesssim \| \A; L_T^\infty H^\sigma \|^{s-1}
\|\A; L^2_T L^\infty\|^{2-s}\]
by the H\"older inequality for the time variable.
\QED
 
\begin{lemma}\label{lem:H2est}\samepage
Let $\sigma>1$ and let $\A\in \M_T^{1,\sigma}\cap L^2_T L^\infty$ satisfy $\Div \A=0$\textup{.}
Let $u \in L_T^\infty H^1$\textup{.}
Then for any $v_0 \in H^2$\textup{,}  there exists a unique weak solution $v$ to 
\eqref{eq:LS}-\eqref{eq:LSinitial}\textup{,}  
and the solution $v$ satisfies the following estimate\textup{:}
\begin{align}\label{eq:H2est}
\| v;L_T^\infty H^2\| &\le C \| v_0;H^2\| \langle \| \A;L_T^\infty \dot{H}^1\| \rangle^4 \nonumber \\
&\quad\times\exp \{ CT^{1/2}\langle T \rangle^l 
\langle  \| \A;\M^{1,\sigma}_T \cap L^2_TL^\infty\| \vee\| u;L^\infty_T H^1\|\rangle^l\}.
\end{align}
Here $l$ is a positive number\textup{.}
Moreover\textup{,} if $u\in C_T H^1$\textup{,} $v$ is a unique strong solution to 
\eqref{eq:LS}-\eqref{eq:LSinitial}\textup{.}
\end{lemma}
\Proof
From the equation \eqref{eq:LS}, we can immediately show the conservation law  
of the $L^2$-norm $\| v(t)\|_2=\|v_0\|_2$.
The uniqueness of the solution clearly follows from this identity.
The existence of the weak solution follows from the a priori estimate \eqref{eq:H2est}.
Indeed the unique existence of the strong solution has already known
 for sufficiently smooth $u$,  $\A$ and $v_0$ (see for example~\cite{K70,K73}).
 Therefore we approximate these functions by a sequence of smooth ones 
 and consider the corresponding sequence of solutions.
 If we extract a star-weakly converging subsequence, then the star-weak limit is a weak solution to 
 \eqref{eq:LS}-\eqref{eq:LSinitial}.
  Therefore we formally prove  \eqref{eq:H2est}.
 Taking Lemma \ref{lem:cov} into account, we estimate 
\[ \| v;H^2_\A\|\equiv\| \DeltaA v \|_2+\langle R \rangle^4 \|v\|_2 \] 
instead of  $\| v ;H^2\|$, where $R\equiv \| \A;L_T^\infty \Dot{H}^1\|$.
Taking the time derivative of $\DeltaA v$ and using the equation \eqref{eq:LS},
we find the equation for $\DeltaA v$:
\begin{equation}\label{eq:DeltaAv}
i\partial_t \DeltaA v = (-\DeltaA +\phi)\DeltaA v+2 \partial_t \A\cdot\nablaA v+[\DeltaA,\phi] v.
\end{equation}
Therefore standard energy method shows that
\[
\|  v ;\revised{L^\infty (t_0,t;H^2_\A)} \|\le \| v_0;H^2_{\A_0} \|
+\| 2 \partial_t \A\cdot\nablaA v+[\DeltaA,\phi] v;\revised{L^1(t_0,t;L^2)}\|,
\]
where $\A_0\equiv \A(t_0)$.
By the Sobolev inequality, we have 
$\| \nabla v \|_r \lesssim \| v;H^2\|^{4\delta} \| v;H^{3/2}_6\|^{1-4\delta}$.
Applying this inequality together with the Young inequality and Lemma  \ref{lem:cov}, we obtain
\begin{align}
&2\| \partial_t \A\cdot\nablaA v \|_2 
\le \| \partial_t \A\|_q \|\nabla v \|_r +\|\partial_t\A\|_2\|\A\|_\infty\|v\|_\infty \nonumber \\ 
&\quad\le\epsilon \| v; H^{3/2}_6\| 
+C\{ \epsilon^{(4\delta-1)/4\delta} \| \partial_t\A;H^{\sigma-1} \|^{1/4\delta}  
+\| \partial_t\A\|_2 \| \A\|_\infty\} \| v;H^2_\A\|.
\end{align}
Here $q,r$ are defined by \eqref{eq:dqr} and
 $\epsilon$ is a positive number which will be determined later.
We can easily handle the term $[\DeltaA,\phi] v$ 
by Lemma \ref{lem:phiest} or the Hardy-Littlewood-Sobolev inequality;
we obtain
\[ \| [\DeltaA,\phi] v \|_2 \lesssim \langle \|\A;\Dot{H}^1\| \rangle \| u;H^1\|^2 \| v;H^2\|. \]
Therefore
\begin{align}
&\| v ; \revised{L^\infty (t_0,t;H^2_\A)} \|  \nonumber \\
&\quad
\le \| v_0;H^2_{\A_0}\|  
+\epsilon  (\revised{t}-t_0)^{1/2} \| v; \revised{L^2 (t_0,t;H^{3/2}_6)}\| \nonumber \\
&\qquad+C\int_{t_0}^{\revised{t}}
\{\epsilon^{(4\delta-1)/4\delta}  \| \partial_t\A;H^{\sigma-1} \|^{1/4\delta}  
+\| \partial_t\A\|_2 \| \A\|_\infty +\langle \|\A;\dot{H}^1\|\rangle\| u;H^1\|^2\} \| v;H^2_\A\|  dt'.
\label{eq:intineq}
\end{align}
Taking  Lemma~\ref{lem:smoothing} into account, we choose the positive number $\epsilon$ so small
that the second term in the right-hand side is absorbed in the left-hand side.
To this end we choose $\epsilon$ such that 
$\epsilon \langle T \rangle^{m +1/2}
\langle \| \A;L_T^\infty H^\sigma \cap L^2_T L^\infty \| 
\vee \|u ; L_T^\infty H^1\| \rangle^m\ll 1$ with $m$ stated in Lemma \ref{lem:smoothing}.
\revised{%
Thus we obtain 
\[ \| v;H^2\|\le C \langle R \rangle^4 \| u_0 ;H^2\|
+C\langle T \rangle^l \langle \| \A;\M^{1,\sigma}_T \| \vee \| u;L^\infty_T H^1\| \rangle^l 
\int_{t_0}^t \langle \| \A \|_\infty \rangle \| v;H^2\| dt',\]
where $l$ is some positive number.}
Applying the Gronwall inequality we obtain \eqref{eq:H2est}.
We proceed to the latter part of the lemma.
\revised{
The weak continuity of $\DeltaA v$ follows from the construction of $v$. 
We shall prove the strong continuity.
We 
take a supreme limit of the both sides of \eqref{eq:intineq} 
as $t\downarrow t_0$ and obtain}
$\limsup_{t \downarrow t_0}\| \DeltaA v(t) \|_2 \le\| \Delta_{A_0} v_0 \|_2$.
Therefore $\textrm{s-}\lim_{t \downarrow t_0}\DeltaA v =\Delta_{A_0} v_0$ in $L^2$.
This argument shows that $\DeltaA v\in C_T L^2$; 
on the other hand, the fact that $v \in C_TL^2$ is similarly proved by  the conservation law of the $L^2$-norm.
Therefore we obtain $v\in C_T H^2$ from Lemma \ref{lem:cov} 
and the fact $\A \in \M^{1,\sigma}_T \subset C_T H^1$.
If we also assume $u \in C_T H^1$, we can show $\partial_t v \in C_T L^2$ 
taking the equation \eqref{eq:LS} into account.
\QED
\begin{definition}
We define the two parameter family of operator $\{ U(t,\tau)\}_{0\le t,\tau \le T}$ 
as the evolution operator of \eqref{eq:LS}.
Namely, $U(t,t_0) v_0$ solves \eqref{eq:LS}-\eqref{eq:LSinitial}.
By the following lemmas, $\{ U(t,\tau)\}$ can be extended as a family in $H^s$ with $s\ge -2$.
We put $K_s\equiv \sup_{0\le t,\tau\le T}\| U(t,\tau);H^s\to H^s\|$.
\end{definition}
\begin{remark}
In view of Lemma \ref{lem:H2est}, we have the estimate
\[ K_2\le C\langle \| \A;L_T^\infty \dot{H}^1\| \rangle^4 
\exp \{ CT^{1/2}\langle T \rangle^l \langle  \| \A;\M^{1,\sigma}_T \cap L^2_TL^\infty\| \vee
\| u;L^\infty_T H^1\|\rangle^l\} \]
with some positive numbers $C,l$.
\end{remark}
\begin{lemma} \label{lem:propertyU}\samepage
Let $-2\le s\le 2$\textup{,} $\sigma>1$ and let 
$\A\in \M_T^{1,\sigma}\cap L^2_T L^\infty$ satisfy $\Div \A=0$\textup{.}
Let $u \in L_T^\infty H^1$\textup{.}
Then we have the following\textup{:}  

\textup{(i)}
$\{ U(t,\tau)\}$ can be uniquely extended to a 
family of operators in $H^s$  
which solves \eqref{eq:LS}\textup{,}
namely $U(t,t_0)v_0$ is a weak $H^s$ solution to \eqref{eq:LS}-\eqref{eq:LSinitial},
and moreover if $s\ge 0$, then $U(t,t_0)v_0$ is a unique weak $H^s$ solution\textup{;}

\textup{(ii)}
if $u \in C_T H^1$\textup{,} $U(t,t_0)v_0$ is a strong $H^s$ solution\textup{;}
 
 \textup{(iii)}
$1\le K_s \le K_2^{|s|/2}$\textup{;}

\textup{(iv)}
if $w$ is a weak $L^2$ solution to the equation
$i\partial_t w=(-\DeltaA +\phi) w +f$\textup{,}
where $f\in L^1_T H^{-2}$\textup{,} then $w$ satisfies the integral form of this equation\textup{,} namely
\[ w(t)=U(t,t_0) w(t_0)-i \int_{t_0}^t U(t,\tau) f(\tau)d\tau. \]
\end{lemma}
 \Proof 
 This lemma can be proved in the same way as in Lemmas 3.2 and 3.3 in \cite{NW05}. 
 \QED
 
\begin{lemma}\label{lem:Hsest}\samepage
Let $(s,\sigma)\in \mathcal{R}_*$ with $2<s<4$ and $\sigma>1$\textup{.}
Let $\A \in \M^{2,\sigma}_T \cap L^2_T L^\infty$ satisfy $\Div \A=0$\textup{.}
Let $u\in \Sigma^{1,s-1}_T$\textup{.}
Then for any $v_0\in H^s$\textup{,} 
there exists a unique weak solution $v$ to \eqref{eq:LS}-\eqref{eq:LSinitial}\textup{,}
and the solution $v$ satisfies the following estimate\textup{:}
\begin{align}
\| v;L^\infty_T H^s\| &\le C K_2 \| v_0;H^s\| 
\langle \| \A;\M^{1,\sigma}_T\|\vee \|u;L^\infty_T H^1\|  \rangle^l \nonumber \\
&\qquad \times\exp \{ C \langle K_2 \rangle^l \langle T \rangle^l  
\langle \|\A ;\M^{2,\sigma}_T \cap L^2_T L^\infty\| \vee \| u;\Sigma^{1,s-1}_T \|  \rangle^l \}. \label{eq:Hsest}
\end{align}
Here $l$ is a positive number\textup{.}
Moreover\textup{,} if $\A \in C_T H^\sigma \cap C^1_T H^{\sigma-1}$ and if 
$u\in C_T H^{s-1}$\textup{,} then $v$ is a strong solution to 
\eqref{eq:LS}-\eqref{eq:LSinitial}\textup{.}
\end{lemma}
\Proof
As we have mentioned  in the proof of Lemma~\ref{lem:H2est}, the existence of a weak solution follows from
the estimate \eqref{eq:Hsest}, and the uniqueness has already been proved in Lemma~\ref{lem:H2est}.
Instead of \eqref{eq:Hsest} itself,  we first assume $5/2< s<4$ and prove the estimate 
 \begin{align}
\| v;L^\infty_T H^s\| &\le C K_2 \| v_0;H^s\| 
\langle \| \A;\M^{1,\sigma}_T\|\vee \|u;L^\infty_T H^1\|  \rangle^l \nonumber \\
&\qquad \times\exp \{ C  \langle K_2 \rangle^l  \langle T \rangle^l  
\langle \|\A ;\M^{2,\sigma}_T\cap L^2_T L^\infty\| \vee 
 \| u;\Sigma^{1,s-3/2}_T \| \rangle^l \}, \label{eq:mHsest}
\end{align}
which is slightly different from but stronger than \eqref{eq:Hsest}.
Then \eqref{eq:Hsest} follows directly from \eqref{eq:mHsest} if $5/2< s<4$  and from 
interpolating \eqref{eq:H2est} and \eqref{eq:mHsest}  if $2<s\le 5/2$. 
Therefore we  assume $5/2<s<4$ and prove \eqref{eq:mHsest},
with dividing the proof into several steps.
In the proof we estimate $\| \partial_t \DeltaA v;H^{s-4}\|$ instead of $\| v;H^s\|$ 
taking the equivalence of these norms  into account. 
To this end we first prove the equivalence of these norms and prepare some inequalities in Step 1.
We also use the smoothing property of the Schr\"odinger equation;
practically, in Step 2 we derive estimates which are  consequences of Lemma \ref{lem:smoothing}.
In Step 3 we apply the estimates obtained in the preceding steps to the equation for $\partial_t \DeltaA v$
and derive an integral inequality for $\| v;H^s\|$, from which the desired estimate 
\eqref{eq:mHsest} follows.
In Step 4 we prove the continuity of the solution, namely the latter part of the lemma.
\Step{1}
We can obtain the following estimates for the solution to \eqref{eq:LS}:
\begin{align}
&\| \partial_t v;H^{s-2}\|+\langle N_1\rangle^\alpha \| v\|_2
\simeq
\| v;H^s\|+\langle N_1 \rangle^\alpha \| v\|_2, \label{eq:equiv1}
\\
&\| \partial_t \DeltaA v;H^{s-4}\|+\langle N_2   \rangle^\alpha \| v\|_2
\simeq
\| v;H^s\|+\langle N_2  \rangle^\alpha \| v\|_2, \label{eq:equiv2}
\end{align}
where
$N_1=\|\A;L^\infty_T H^\sigma \|\vee \|u; L^\infty_T H^1\|$,
$N_2=\|\A;\M^{1,\sigma}_T \|\vee \|u; L^\infty_T H^1\|$.
To obtain \eqref{eq:equiv1} and \eqref{eq:equiv2}, 
we have only to use  the following inequalities together with Lemma \ref{lem:cov} 
and standard interpolation inequality:
\begin{gather} 
\| \nablaA  v;H^{s-2}\| \lesssim \langle \| \A; H^\sigma \| \rangle \| v;H^{s-1}\|, \label{eq:nablaAv} \\
\| \phi v;H^{s-2}\| \lesssim \| u; H^1\|^2 \| v;H^{s-1}\|, \label{eq:phivs-2} \\
\| V_2(\A,\partial_t \A,v);H^{s-4}\| \lesssim 
\langle N_2 \rangle^3
\| v;H^{s-1}\|, \label{eq:V2}
\end{gather}
where $V_2(\A, \B,v)=P\B \cdot \nablaA  v$. 
We can prove these inequalities by the Sobolev inequality and Lemma \ref{lem:phiest}.
We remark that \eqref{eq:nablaAv} and \eqref{eq:phivs-2} holds valid for $0\le s<4$ (for the proof we use duality argument) and that $V_2$ is a continuous mapping from 
$H^\sigma \times H^{\sigma-1}\times H^{s-1}$ to $H^{s-4}$.

\Step{2}
We prove
\begin{equation}\label{eq:smoothing2}
\| \DeltaA v; L^2_T H^{s-5/2}_6 \|\le
C \langle T \rangle^m  \langle N_3 \rangle^m \| v;L^\infty_T H^s \|,
\end{equation}
where $N_3=\|\A;\M^{1,\sigma}_T \cap L^2_T L^\infty  \|\vee \|u; L^\infty_T H^{1\vee (s-5/2)}\|$
and $m$ is a positive number.
We remark that the constant $C$ does not depend on $T$ as well as $v$ and $N_3$.
Applying Lemma \ref{lem:smoothing} to \eqref{eq:DeltaAv}, we have
\begin{align}
\| \DeltaA v; L^2_T H^{s-5/2}_6 \| &\lesssim
\langle T \rangle^m  \langle N_3 \rangle^m  \| \DeltaA v;L^\infty_T H^{s-2} \| \nonumber \\
&\quad+T^{1/2} \| 2 \partial_t \A\cdot\nablaA v+[\DeltaA,\phi] v;L^2_T H^{s-3}\|.
\end{align}
Therefore we obtain \eqref{eq:smoothing2} by \eqref{eq:cov} if we prove the following estimates:
\begin{align}
\| \partial_t \A\cdot\nablaA v;H^{s-3}\| 
&\lesssim
\| \partial_t \A;H^{\sigma-1} \| \langle \| \A; H^\sigma \| \rangle \| v;H^s\|,  \label{eq:Anablav}
\\
\| |u|^2 v ;H^{s-3}\| 
&\lesssim
\| u; H^{1\vee (s-5/2)}\|\|u;H^1\| \| v;H^s\|, \\
\| \nabla \phi \cdot \nablaA v; H^{s-3} \| 
&\lesssim
\| u; H^1\|^2 \langle \| \A; H^\sigma \| \rangle \| v;H^s\|. 
\end{align}
These estimates can be proved by the Sobolev inequality, together with the Leibniz rule if $s\ge 3$, and moreover the duality argument if $s<3$.
We need the condition $\sigma\ge (2s-1)/4$ to bound $\|\partial_t\A\cdot\A v;H^{s-3}\|$ by the right hand side of \eqref{eq:Anablav}.
Replacing $\DeltaA v$ in \eqref{eq:smoothing2} by $-i\partial_t v +\phi v$ and using Lemma \ref{lem:phiest},
we also obtain 
\begin{equation}\label{eq:smoothing3}
\| \partial_t v; L^2_T H^{s-5/2}_6 \|\le
C \langle T \rangle^m  \langle N_3 \rangle^m  \| v;L^\infty_T H^s \|.
\end{equation}

\Step{3}
We take the time derivative of \eqref{eq:DeltaAv}. Then we obtain the following Schr\"odinger equation for
$\partial_t\DeltaA v$:
\begin{align}
i\partial_t^2  \DeltaA v &= (-\DeltaA +\phi)\partial_t\DeltaA v \nonumber \\
&\quad+2 \partial_t \A\cdot\nablaA \DeltaA v 
+\partial_t \phi \cdot \DeltaA v 
+2\partial_t^2  \A\cdot\nablaA  v 
-2i |\partial_t \A |^2 v 
+2 \partial_t \A\cdot\nablaA \partial_t v \nonumber \\ 
&\quad-2 \Re (\partial_t u \cdot \Bar{u}) v 
+2 \nabla \partial_t \phi \cdot \nablaA v 
-2i \nabla \phi \cdot \partial_t \A \cdot v 
-|u|^2 \partial_t  v 
+2 \nabla \phi \cdot \nablaA \partial_t v  \nonumber \\ 
&\equiv (-\DeltaA +\phi)\partial_t\DeltaA v + \sum_{j=1}^{10} f_j, \label{eq:DtDeltaAv}
\end{align}
where we have used the relation 
$\partial_t\nablaA=\nablaA\partial_t-i\partial_t\A$,
$\partial_t\DeltaA=\DeltaA\partial_t-2i\partial_t\A\cdot\nablaA$,
$\Delta \partial_t\phi=-2\mbox{Re} \partial_tu\cdot \Bar{u}$, and 
$[\DeltaA,\phi]v=\Delta\phi\cdot v+2\nabla \phi\cdot\nablaA v$.
By the Duhamel principle or Lemma \ref{lem:propertyU}, 
we rewrite this equation into integral form as
\begin{equation}
\partial_t\DeltaA v = U(t,t_0)[\partial_t\DeltaA v]_{t=t_0} 
-i\int_{t_0}^t U(t,\tau) \sum_{j=1}^{10} f_j(\tau) d\tau. \label{eq:intDtDeltaAv}
\end{equation}

Therefore we have
\begin{equation} \label{eq:W}
W(t) \lesssim K_2 \Bigl\{ W(t_0)+
\int_{t_0}^t \sum_{j=1}^{10}\| f_j(\tau) ; H^{s-4}\| d\tau \Bigr\},
\end{equation}
where $W(t) =\| v(t);H^s \|+\langle N_2 \rangle^\alpha \| v(t)\|_2$.
We have used \eqref{eq:equiv2}, Lemma \ref{lem:propertyU}, 
 and the conservation law of $\| v\|_2$.
We estimate $f_j$, $j=1,\dots, 10$, and obtain the following estimates:
\begin{align}
\| f_1;H^{s-4} \| &\lesssim \| \partial_t \A; H^{\sigma-1} \|  \langle \| \A; H^\sigma \| \rangle
\| \DeltaA v; H^{s-5/2-\delta}_6 \|, 
\label{eq:f1} \\
\| f_2;H^{s-4} \| &\lesssim \| \partial_t u; H^{-1}\| \| u;H^1\|  \| \DeltaA v; H^{s-2} \|, 
\label{eq:f2} \displaybreak[0] \\ 
\| f_3;H^{s-4} \| &\lesssim  
\| \partial_t^2 \A; H^{\sigma-2} \|  \langle \| \A; H^\sigma \| \rangle \| v; H^s\|, 
\label{eq:f3} \\
\| f_4;H^{s-4} \| &\lesssim \| \partial_t \A; H^{\sigma-1} \|^2 \| v; H^s\|, 
\label{eq:f4} \displaybreak[0]  \\
\| f_5;H^{s-4} \| &\lesssim \| \partial_t \A; H^{\sigma-1} \| \langle \| \A; H^\sigma \| \rangle 
\| \partial_t v; H^{s-5/2-\delta}_6 \|,  
\label{eq:f5}  \\
\| f_6;H^{s-4} \| &\lesssim \| \partial_t u; H^{s-7/2}\| \| u; H^{s-3/2} \|  \| v; H^s\|,  
\label{eq:f6}  \displaybreak[0] \\
\| f_7;H^{s-4} \| &\lesssim \| \partial_t u; H^{s-7/2}\| \| u;H^{s-3/2}\| 
\langle \| \A; H^\sigma \| \rangle \| v; H^s \|, 
\label{eq:f7}  \\
\| f_8;H^{s-4} \| &\lesssim \| \partial_t \A \|_2 \| u;H^1\|^2 \| v; H^s\|, 
\label{eq:f8}  \displaybreak[0] \\
\| f_9;H^{s-4} \| &\lesssim \| u; H^1 \|^2 \| \partial_t v; H^{s-2}\|,  
 \label{eq:f9}  \\
\| f_{10};H^{s-4} \| &\lesssim \langle \| \A; H^\sigma \| \rangle \| u;H^1\|^2 \| \partial_t v; H^{s-2}\|.
\label{eq:f10} 
\end{align}
The inequality \eqref{eq:f1} is obtained by the use of 
the estimates 
\revised{%
\[ \|\partial_t \A\cdot (\nablaA w);H^{s-4}\| =\|\nablaA \cdot (\partial_t \A w);H^{s-4}\| 
\lesssim \langle \|\A;H^1\|\rangle\|\partial_t \A w;H^{s-3}\| \]
for any $w\in H^{s-3}$, and }   
\[ \| \partial_t \A \DeltaA v;H^{s-3}\| \lesssim 
\| \partial_t \A ;H^{(s-3)\vee 2\delta}\|  \| \DeltaA v; H^{s-5/2-\delta}_6 \|. \]
\revised{The second} inequality is obtained by the Leibniz rule and the Sobolev inequality for $s=3+2\delta$ and for $s=4$, 
together with duality argument for $s=5/2$, and it is generalized by interpolation for $5/2<s<4$.
The inequality \eqref{eq:f5} is obtained in the same way. 
We can obtain the other inequalities principally by the H\"older, the Leibniz, and the Sobolev inequalities,
together with duality argument if Sobolev spaces of negative order appear.
We also use Lemma \ref{lem:phiest} for the proof of \eqref{eq:f2}, \eqref{eq:f7}, \eqref{eq:f8}
and  \eqref{eq:f10}.
The condition $\sigma\ge(2s-1)/4$ is needed to bound the terms $f_3$ and $f_4$.
We substitute these estimates into \eqref{eq:W}, and apply the Gagliardo-Nirenberg inequality 
$\| \DeltaA v; H^{s-5/2-\delta}_6 \| \lesssim 
\|\DeltaA v; H^{s-2} \|^{2\delta} \| \DeltaA v; H^{s-5/2}_6 \|^{1-2\delta}$
and the corresponding inequality for $\partial_t v$ together with the Young inequality.
Then we obtain
\begin{align*}
W(t) &\le K_2 \Bigl\{ W(t_0)+\epsilon \int_{t_0}^t \bigl(
\| \DeltaA v; H^{s-5/2}_6 \|+\| \partial_t v; H^{s-5/2}_6 \|\bigl) d\tau \\
&\quad+C\epsilon^{-2\delta/(1-2\delta)} \int_{t_0}^t \langle N_4 \rangle^l W(\tau) d\tau \Bigr\} \\
&\le K_2 W(t_0) +C \varepsilon K_2 \langle T \rangle^{m+1/2}  \langle N_3 \rangle^m  \| v;L^\infty(t_0,t; H^s) \| \\
&\quad+C K_2 \epsilon^{-2\delta/(1-2\delta)} \int_{t_0}^t \langle N_4 \rangle^l W(\tau) d\tau,
\end{align*}
where $N_4=\| \A;\M^{2,\sigma}_T\cap L^2_TL^\infty \| \vee \| u;\Sigma^{1,s-3/2}_T \|\ge N_3$ 
and $\epsilon$ and $l$ are positive numbers. We choose $\epsilon$ such that
$C \varepsilon K_2 \langle T \rangle^{m+1/2}  \langle N_3 \rangle^m \le 1/2$;
with this choice we obtain
\[ W(t)\le 2 K_2 W(t_0) 
+C \langle T \rangle^l \langle N_4 \vee K_2\rangle^l  
\int_{t_0}^t  W(\tau) d\tau. \]
We note that the value of $l$ may differ from that in the previous estimate.
Applying the Gronwall lemma, we obtain 
$W(t)\le 2 K_2 W(t_0) \exp \{ C \langle T \rangle^l \langle N_4 \vee K_2\rangle^l \}$,
from which we can conclude \eqref{eq:mHsest}.

\Step{4}
We prove the continuity of $v$ under the additional assumptions 
$\A\in C_T H^\sigma \cap C^1_T H^{\sigma-1}$ and $u\in C_T H^{s-1}$.
We first remark that $v\in C_T H^{s-1}$ since $v\in C_T H^2 \cap L^\infty_T H^s$ 
by virtue of Lemma~\ref{lem:H2est} and \eqref{eq:Hsest}. 
By the assumption, $U(t,\tau)$ is $H^{s-4}$-strongly continuous with respect to the parameters $t,\tau$.
Moreover, $\sum_{j=1}^{10} f_j \in L^1 H^{s-4}$ by Step 3.
Applying the Lebesgue convergence theorem to \eqref{eq:intDtDeltaAv},
we can prove $\partial_t \DeltaA v \in C_T H^{s-4}$.
Then we find that $\DeltaA \partial_t v =\partial_t \DeltaA v+2i\partial_t\A\cdot\nablaA v\in C_T H^{s-4}$
taking \eqref{eq:V2} into account.
Therefore it follows from Lemma~\ref{lem:cov} that $\partial_t v \in C_T H^{s-2}$. 
We go back to the equation \eqref{eq:LS} and conclude $v \in C_T H^s$ 
by Lemmas~\ref{lem:phiest} and \ref{lem:cov}.
\QED


\section{Unique existence of local solutions}\label{sec:local}
In this section we uniquely solve (MS-C) time locally.
To this end we consider the following linearized equation:
\begin{alignat}{2}
& i\partial_t v=(-\DeltaA+\phi(u))v, & \quad & v(0)=u_0,  \label{eq:v} \\
& (\square+1)\B=P\J(u,\A)+\A, & \quad & B(0)=\A_0,\partial_t\B(0)=\A_1,\label{eq:B}
\end{alignat}
where $\phi(u)=(-\Delta)^{-1}|u|^2$ and $\J(u,\A)=2\Im \Bar u \nablaA u$,
with the assumptions $\Div \A=0$ and $(u_0,\A_0,\A_1)\in X^{s,\sigma}$.
We often consider the equations with $(u,\A, v,\B)$ replaced by $(u',\A',v',\B')$ 
and $(u_0,\A_0,\A_1)$ by $(u_0',\A_0',\A_1')$.
In such a case we often abbreviate $\phi(u'), \J(u',\A')$ to $\phi', \J'$.
If we define the mapping 
\[ \Phi: (u,\A)\mapsto (v,\B),\]
then the fixed points of $\Phi$ solve (MS-C).

\begin{proposition}\label{prop:localexistence} \samepage
Let $11/8\le s \le 2$ and $1<\sigma\le \min\{3s/2; 2s-3/4\}$ with $(s,\sigma)\neq (2,3)$\textup{.}
Then for any $(u_0,\A_0,\A_1)\in X^{s,\sigma}$\textup{,} there exists $T>0$ such that 
\textup{(MS-C)} with \eqref{eq:IDC} has a unique solution satisfying 
$(u,\A, \partial_t \A)\in C_T X^{s,\sigma}$ and
$(u,\A)\in L^2_T (H^{s-1/2}_6 \oplus L^\infty)$\textup{.}
Moreover\textup{,} the total energy 
\[ \mathcal{E}=\| \nablaA u\|_2^2+\frac12 \{ \|\nabla \phi\|_2^2 +\| \nabla \A \|_2^2+\| \partial_t \A \|_2^2\} \]
does not depend on $t$\textup{.}
\end{proposition}
\Proof
Let $0<T<1$ and $R_1,R_2, R_3>1$.
We define the metric space $\mathcal{B}$ with metric induced from the norm 
$\|\cdot;\mathcal{B}\|$ as 
\begin{align*}
\mathcal{B}&=\{(u,\A); u\in L^\infty_T H^s \cap L^2_T H^{s-1/2}_6,
 \A \in \M^{1,\sigma}_T \cap L^2_T L^\infty, \Div \A =0, \\
&\qquad\quad \| u;L^\infty_T H^s \|\le R_1,
\| u; L^2_T H^{s-1/2}_6 \| \le R_2,
\| \A;\M^{1,\sigma}_T\cap L^2_T L^\infty \|\le R_3 \}, 
\end{align*}
\[ \| (u,\A);\mathcal{B}\|= \| u;L^\infty_T L^2\| \vee \| \A; L^\infty_T H^{1/2}\cap L^4_T L^4\|. \]
We can easily show that this metric space is complete.
We shall show that $\Phi$ is a contraction mapping defined on $\mathcal{B}$.
On account of Lemmas \ref{lem:StKG}, \ref{lem:nlmaxwell}, 
\ref{lem:smoothing} and \ref{lem:propertyU}, 
we have the following estimates for $(u,\A) \in \mathcal{B}$ and $(v,\B)=\Phi(u,\A)$:
\begin{gather*}
\| v;L^\infty_T H^s \|\le C \| u_0;H^s\| R_3^{2s} \exp \{C T^{1/2}  (R_1 \vee R_3)^l\}, \\
\| v;L^2_T H^{s-1/2}_6 \| \le C(R_1 \vee R_3)^m \| v;L^\infty_T H^s \|, \\
\| \B;\M^{1,\sigma}_T \cap L^2_T L^\infty\| 
\le C\| (\A_0,\A_1);H^\sigma\oplus H^{\sigma-1}\| + CT^{1/4} R_3 (R_1 \vee R_2)^2.
\end{gather*}
To prove the last inequality, we have also used the Sobolev type embedding 
$H^{\sigma-2/q}_r \embedding L^\infty$, where $q, r$ are defined in \eqref{eq:dqr}.
Therefore we can show that $\Phi$ is a mapping from $\mathcal{B}$ to itself
if we choose $R_1,R_2, R_3$ such that
\[ R_3 \ge 2C\| (\A_0,\A_1);H^\sigma\oplus H^{\sigma-1}\|, \quad
R_1 \ge 2C  \| u_0;H^s\| R_3^{2s}, \quad
R_2 \ge C(R_1 \vee R_3)^m R_1\]
and if we choose $T$ such that
\[ \exp \{C T^{1/2}  (R_1 \vee R_3)^l\}\le 2,  \quad CT^{1/4}  (R_1 \vee R_2)^2 \le 1/2. \]
Next we estimate the difference of $(v,\B)=\Phi(u,\A)$ and $(v',\B')=\Phi(u',\A')$.
Taking the difference of the equations for $v$ and $v'$, we obtain
\begin{align*}
i\partial_t(v-v')&=(-\DeltaA +\phi) (v-v') +2i (\A-\A')\cdot \nabla v'  \\
&\quad+ (\A-\A')\cdot (\A+\A') v' +(\phi -\phi') v'.
\end{align*}
By usual $L^2$-estimate together with the Sobolev inequality,
Lemma \ref{lem:phiest} and the H\"older inequality both for time and space variables, 
we obtain
\begin{align*}
&\| v-v';L^\infty_T L^2\|  \nonumber \\
&\quad\le \| 2i (\A-\A')\cdot \nabla v'  
+ (\A-\A')\cdot (\A+\A') v' +(\phi -\phi') v'; L^1_T L^2 \| \nonumber \\
&\quad \lesssim \| \A-\A'; L^4_T L^4 \| \{ T^{3/8}\| \nabla v'; L^{8/3}_T L^4 \|
+T^{3/4}\| \A+\A';L^\infty L^6\| \| v';L^\infty_T L^{12} \| \} \nonumber \\
&\qquad+T \| u-u';L^\infty_T L^2\|  \| u+u';L^\infty_T H^1 \| \| v';L^\infty_T H^1 \| \nonumber \\
&\quad\lesssim  T^{3/8} (R_1\vee R_2 \vee R_3)^2 
\| (u-u', \A-\A'); \mathcal{B}\|. \nonumber 
\end{align*}
Here we have used the interpolation
$L^{8/3}_T H^{s-3/8}_4=(L^\infty_T H^s, L^2_T H^{s-1/2}_6)_{[3/4]}$.
We proceed to 
the estimate of $\B-\B'$.
We apply Lemma \ref{lem:StKG} to the difference of the equations for $\B$ and $\B'$
taking the relation $P\nabla=0$ into account. Then
\begin{align*}
&\| \B-\B';L^\infty_T H^{1/2}\cap L^4_T L^4 \|  \\
&\quad\lesssim \| \A-\A'; L^1_T H^{-1/2}\|
+\| P(\J-\J'); L^{4/3}_T L^{4/3} \| \\
&\quad\lesssim 
T\| \A-\A'; L^\infty_T H^{1/2}\|
+T^{3/8} \| \nabla (u+u'); L^{8/3}_T L^4 \| \| u-u';L^\infty_T L^2 \| \\
&\qquad+T\| \A; L^\infty_T L^6\| \| u-u';L^\infty_T L^2 \| \| u+u';L^\infty_T L^{12} \|  \\
&\qquad+T^{1/2} \| \A-\A'; L^4_T L^4 \| \| u';L^\infty_T L^4 \|^2 \\
&\quad \lesssim T^{3/8} (R_1\vee R_2\vee R_3)^2 \| (u-u', \A-\A'); \mathcal{B}\|,
\end{align*}
where we have used \revised{the same interpolation relation as above.}
Therefore we obtain
\begin{equation}\label{eq:contraction}
\| (v-v', \B-\B'); \mathcal{B}\|\le (1/2)\| (u-u', \A-\A'); \mathcal{B}\|
\end{equation}
for sufficiently small $T>0$.
Therefore, $\Phi$ is a contraction mapping with the choice of $T, R_1,R_2,R_3$ mentioned above,
from which we conclude the unique existence of the solution.
Moreover $(u,\A,\partial_t \A)\in C_T X^{s,\sigma}$. 
Indeed,  
$(\A,\partial_t \A)\in C_T (H^\sigma \oplus H^{\sigma-1})$ by Lemma \ref{lem:StKG},
and $u \in C_T H^1$ since solutions to \eqref{eq:v} belong to $\Sigma^{1,s}_T \subset C_T H^1$
by Lemma~\ref{lem:H2est};
therefore going back to the Schr\"odinger part, 
we obtain $\revised{u}\in C_T H^s \cap C^1_T H^{s-2}$ by Lemma \ref{lem:propertyU}.
Finally, we prove the conservation of total energy.
For $H^2$-strong solutions, this follows from direct computation.
For a solution $(u, \A)$ with lower regularity, 
we consider a sequence of $H^2$-solutions $\{(u^j, \A^j)\}_j$ which is an approximation to $(u, \A)$.
As we obtained \eqref{eq:contraction}, we can prove that $\{(u^j, \A^j,\partial_t \A^j)\}_j$
converges to $(u,\A , \partial_t \A)$ in $L^\infty_T X^{0,1/2}$.
Since $\{(u^j, \A^j,\partial_t \A^j)\}_j$ is bounded in $L^\infty_T X^{s,\sigma}$, this sequence
actually converges to $(u,\A , \partial_t \A)$ in $L^\infty_T X^{1,1}$.
Therefore the conservation of total energy holds also for \revised{$(u,\A)$}.
\QED

\begin{proposition}\label{prop:regularity} \samepage
Let $(s,\sigma)\in \mathcal{R}$ with $s>2$ and $\sigma >1$\textup{.}
Then for any $(u_0,\A_0,\A_1)\in X^{s,\sigma}$\textup{,} 
the solution to \textup{(MS-C)} with \eqref{eq:IDC} obtained in Proposition \textup{\ref{prop:localexistence}}
actually belongs to $C_T X^{s,\sigma}$\textup{.}
\end{proposition}
\Proof
Firstly let $(s,\sigma)\in \mathcal{R}$ satisfy $2<s< 4$ and $\sigma<3$.
Then the unique solution obtained by Proposition \ref{prop:localexistence}
belongs to $C_T X^{2,\sigma}$.
Using Lemma \ref{lem:Hsest} at most twice, we can prove that $u\in C_T H^s \cap C^1_T H^{s-2}$
(remark that $\partial_t^2 \A=\Delta\A+P\J\in C_T H^{\sigma-2}$ by virtue of Lemma \ref{lem:commutator}),
and hence we can prove the proposition for such $(s,\sigma)$.
Next we apply Lemmas \ref{lem:StKG}, \ref{lem:nlmaxwell} and \ref{lem:smoothing}  
to the solution and obtain the proposition for $(s,\sigma)\in \mathcal{R}$ with $2<s< 4, \sigma\ge3$.
The proposition for $s\ge 4$ has already been obtained in~\cite{NW05}.
\QED
\thmskip
\Proofof{Theorem \textup{\ref{thm:MS-C}}}
We can prove the existence of solutions by combining Propositions \ref{prop:localexistence}-\ref{prop:regularity}.
The uniqueness without auxiliary conditions $(u,\A)\in L^2_T (H^{s-1/2}\oplus L^\infty)$
is a consequence of Corollary \ref{cor:apriori} in \S \ref{sec:global}.
The continuous dependence of solutions on initial data will be proved by
Proposition \ref{prop:continuous} in \S \ref{sec:data}. \QED


\section{Global existence of solutions}\label{sec:global}
\begin{lemma} \label{lem:apriori}\samepage
Let $11/8\le s\le 2$\textup{,} $1<\sigma\le 10/9$ and $(u,\A,\partial_t \A) \in C_TX^{s,\sigma}$ be a solution to \textup{(MS-C)}
obtained in Proposition \textup{\ref{prop:localexistence}.}
Then the following estimates hold\textup{.} 
\begin{gather}
\| (u, \A, \partial_t\A);L^\infty_T(H^1\oplus \Dot{H}^1\oplus L^2) \| \le C, \label{eq:apriori1} \\
\| \A;L^\infty_TL^2\| \le C\langle T\rangle, \label{eq:apriori2} \\
\| \A;L^q_TL^r\| \le C\langle T\rangle^2, \label{eq:apriori3} \\
\| u;L^2_T H^{1/2-\delta}_6\| \le C\langle T\rangle^3, \label{eq:apriori4}  \\
\| \A;\M_T^{1,\sigma}\cap L^2_T L^\infty \| \le C\langle T\rangle^4. \label{eq:apriori5} 
\end{gather}
Here $q,r$ and $\delta$ are given in \eqref{eq:dqr}\textup{.}
The constants $C$ depend only on $\sigma$ and
$\| (u_0,\A_0,\A_1);X^{\revised{s,\sigma}}\|$\textup{.}\end{lemma}
\Proof
We easily obtain \eqref{eq:apriori1} by the conservation laws of charge and energy,
\eqref{eq:apriori2} by applying \eqref{eq:apriori1} after differentiating and integrating $\A$ with respect to $t$.
Next we apply Lemma~\ref{lem:StKG} to \eqref{eq:MSC2}
and obtain
\begin{align*}
\| \A;L^q_TL^r\| &\lesssim \| (\A_0,\A_1);H^1\oplus L^2\| + \| \A ;L^1_T  H^{2/q-1} \|
+\| P\J ;L^{q'}_T H^{4/q-1}_{r'} \|. 
\end{align*}
The second term in the right-hand side is bounded by $C\langle T \rangle^2$ because of  \eqref{eq:apriori2},
and the third term is bounded by
$T^{1/q'}\| u ;L^\infty_T H^1 \|^2 \langle \| \A;L^\infty_T H^1\| \rangle$
by the use of Lemma \ref{lem:commutator} together with the Leibniz rule and the Sobolev inequality.
Hence this term is also bounded by $C\langle T \rangle^2$ and \eqref{eq:apriori3} has been proved.
In order to obtain \eqref{eq:apriori4}, we apply Lemma~\ref{lem:KT} to \eqref{eq:MSC1}.
Then
\begin{align*}
\| u;L^2_T H^{1/2-\delta}_6\| &\lesssim \| u;L^\infty_T H^1\|
+T^{1/2}\| 2i \A\cdot\nabla u+|\A |^2 u +\phi u;L^2_T H^{-2\delta}\|  \\
&\lesssim\langle T\rangle \| u;L^\infty_T H^1\| 
 \langle \| A;L^q_TL^r\| +\| A;L^\infty_T\Dot{H}^1\|^2+ \| u;L^\infty_T H^1\|^2\rangle.
\end{align*}
The right-hand side is estimated by $C\langle T \rangle^3$ by the previous estimates. 
Therefore  \eqref{eq:apriori4} has been proved.
We go back to the Maxwell part and again apply Lemma~\ref{lem:StKG} to \eqref{eq:MSC2}.
Then
\begin{align*}
&\| \A;\M^{1,\sigma}_T\revised{\cap L^q_T H^{\sigma-2/q}_r} \|  \\
&\quad\lesssim \| (\A_0,\A_1);H^\sigma \oplus H^{\sigma-1}\|
+ T\| \A ;L^\infty_T H^{\sigma-1} \| +\| P\J ;L^{6/5}_T H^{\sigma-2/3}_{3/2} \|. 
\end{align*}
By the assumption, $\sigma-2/3\le 1/2-\delta$.
Therefore the last term in the right-hand side is bounded by 
\[ T^{1/3} \| u; L^2_T H^{1/2-\delta}_6 \| \| u;L^\infty_T H^1 \| \langle \| A;L^\infty_T\Dot{H}^1\| \rangle
\lesssim \langle T \rangle^{3+1/3}. \]
If we use the estimate $\| \A;L^2_T L^\infty \| \lesssim T^{1/r} \| \A ;L^q_T H^{\sigma-2/q}_r\|$,
which is obtained by the Sobolev inequality, we can show \eqref{eq:apriori5}. \QED

\begin{corollary}\label{cor:apriori}\samepage
Let $11/8\le s\le 2$\textup{,} $1<\sigma\le 10/9$ and let $(u,\A)$ be a solution to \textup{(MS-C)}
satisfying $(u,\A,\partial_t\A)\in C_T X^{s,\sigma}$\textup{.}
Then $u\in L^2_T H^{s-1/2}_6$ and $\A \in L^2_T L^\infty$.
\end{corollary}
\Proof
If we check the proof of Lemma \ref{lem:apriori}, 
we find that we can prove $\A \in L^2_T L^\infty$ under the assumption 
that $(u,\A)$ satisfies (MS-C) and that $(u,\A,\partial_t\A)\in C_T X^{1,\sigma}$.
Once we have proved  $\A \in L^2_T L^\infty$, we immediately obtain  $u\in L^2_T H^{s-1/2}_6$
by Lemma~\ref{lem:smoothing}. \QED
\thmskip
\Proofof{Theorem \textup{\ref{thm:global}}}
We first consider the case $s\le 2$ and $\sigma\le 10/9$.
By Lemma \ref{lem:apriori}, $\| u;L^\infty_T H^1\|$ and $\| \A; \M^{1,\sigma}_T \cap L^2_T L^\infty\|$
are finite as long as the solution exists in $0<t<T$.
Therefore $\| u; L^\infty_T H^s\|$ is also finite by virtue of Lemma \ref{lem:propertyU}.
This implies the global existence.
For general case we have only to recover the regularity by using Propositions 
\ref{prop:localexistence} and \ref{prop:regularity}.
\QED
 

\section{Continuous dependence on initial data}\label{sec:data}
In this section we shall complete the proof of Theorem \ref{thm:MS-C}
by proving the continuous dependence of solutions on data.
The argument here is essentially based on Bona-Smith~\cite{BS75}.

\begin{lemma}\label{lem:diffest}\samepage
Let $11/8 \le s<4$\textup{,} $\sigma>1$ and let $(s,\sigma)\in\mathcal{R}$ with 
$(s+1,\sigma)\in\mathcal{R}_*$\textup{.}
Let $(u,\A)$ and $(u',\A')$ be solutions to \textup{(MS-C)} defined on $[0,T]$ 
with the initial data $(u_0,\A_0,\A_1)\in X^{s,\sigma}$ and 
$(u_0',\A_0',\A_1')\in X^{s+1,\sigma}$ respectively\textup{.}
Let $(u,\A)$ satisfy the estimate
$\| u;\Sigma^{1,s}_T \cap L^2_T H^{s-1/2}_6 \| \vee
\| \A;\M^{1,\sigma}_T \cap L^q_T H^{\sigma-2/q}_r  \|  \le R$
and $(u',\A')$ satisfy the same estimate with $(u,\A)$ replaced by $(u',\A')$.
Then we have the following estimates\textup{:} \pagebreak[0]
\begin{align}
&\| u_{-};L^\infty_T H^s \cap L^2_T H^{s-1/2}_6 \| \vee 
\| \A_{-};\M^{1,\sigma}_T \cap L^q_T H^{\sigma-2/q}_r  \| \nonumber \\
&\quad \le
C \| (u_0,\A_0,\A_1)_{-}; X^{s,\sigma} \|
+ C\|\A_{-};L^\infty_T H^1\cap L^q_T H^{\sigma-2/q-\delta}_r  \| \| u';\Sigma^{1,s+1}_T \|, \label{eq:diffest1} \\
&\| \A_{-};\M^{1,\sigma}_T \cap L^q_T H^{\sigma-2/q}_r  \| \nonumber \\
&\quad \le
C \| (\A_0,\A_1)_{-}; H^\sigma \oplus H^{\sigma-1} \|
+ C\| u_{-};L^\infty_T H^s \cap L^2_T H^{s-1/2}_6 \|. \label{eq:diffest2}
\end{align} 
Moreover\textup{,} 
let $s>11/8$, $0<\sigma-1-\delta\le \revised{s-11/8}$\textup{.} 
Then we also have the following estimates\textup{:} \pagebreak[0]
\begin{equation}
\| u_{-};L^\infty_T H^{s-1} \| \vee 
\| \A_{-};\M^{1,\sigma-\delta}_T \cap L^q_T H^{\sigma-2/q-\delta}_r  \| 
\le C \| (u_0,\A_0,\A_1)_{-}; X^{s-1,\sigma-\delta} \|. \label{eq:diffest3}
\end{equation}
Here $u_{-}=u-u'$ etc\textup{.,} $\delta, q, r$ are defined in \eqref{eq:dqr}
and \revised{the} constants $C$ depend on $R,T,s$ and~$\sigma$\textup{.}
\end{lemma}
\Proof
It suffices to show \eqref{eq:diffest1}-\eqref{eq:diffest3} for sufficiently small $T=T(R)$;
if not, we divide the interval $[0,T]$ into small subintervals and repeatedly use the estimates
obtained for short intervals. Hence we may assume $0<T<1$ without loss of generality.
We begin with the estimate of the Schr\"odinger part.
Taking the difference of the equations for $u$ and $u'$, we have
\begin{align}
i\partial_t u_{-}&=(-\DeltaA +\phi) u_{-}+2i\A_{-}\cdot\nabla u'+\A_{-}\cdot \A_+ u'+\phi_{-}u'
\nonumber \\
&\equiv(-\DeltaA +\phi) u_{-} +\sum_{j=1}^3 g_j,\label{eq:diff}
\end{align}
where $\A_+=\A+\A'$.
We also need the time derivative of \eqref{eq:diff}:
\begin{align}
i\partial_t^2 u_{-}&=(-\DeltaA +\phi) \partial_t u_{-} \nonumber \\
&\quad
+2i \partial_t\A\cdot\nabla u_{-} +2\partial_t\A\cdot \A u_{-} +\partial_t \phi u_{-} 
+2i \partial_t\A_{-}\cdot\nabla u'+2i \A_{-}\cdot\nabla \partial_t u' \nonumber \\
&\quad
+\partial_t\A_{-}\cdot \A_+ u' + \A_{-}\cdot \partial_t \A_+ u'  +  \A_{-}\cdot \A_+ \partial_t u' 
+\partial_t \phi_{-} u' + \phi_{-} \partial_t u' \nonumber \\
&\equiv(-\DeltaA +\phi) \partial_t u_{-} +\sum_{j=4}^{13} g_j \label{eq:diff2}
\end{align}
In the following, we estimate $\partial_t u_{-}$ instead of $u_{-}$ itself 
in order to obtain \eqref{eq:diffest1} and \eqref{eq:diffest3}.
To this end we introduce here an inequality which shows the equivalence of norms
$\| u_{-};H^s \|$ and $\| \partial_t u_{-};H^{s-2}\|$.
Namely for $(s,\sigma)\in \mathcal{R}_*$ with $s>1/2$ we have
\begin{equation} \label{eq:equiv3}
\| u_{-};H^s \|+C(R) \{ \| u_{-}\|_2 +\|\A_{-};H^\sigma\|\}
\simeq \| \partial_tu_{-};H^{s-2} \|+C(R) \{ \| u_{-}\|_2 +\|\A_{-};H^\sigma\|\}.
\end{equation}
We can prove this inequality in the same way as we proved \eqref{eq:equiv1},
namely we use a trivial modification of \eqref{eq:V1} together with Lemma \ref{lem:phiest}.
We refer to Lemma~6.1 in~\cite{NW05} for detail.
Next we apply Lemma~\ref{lem:smoothing} to \eqref{eq:diff} and obtain
\begin{equation}\label{eq:diffsmoothing}
\| u_{-};L^2_T H^{s-1/2}_6\| \le C(R) \{ \| u_{-};L^\infty_T H^s\| \vee 
\| \A_-;L^\infty_T H^\sigma \cap L^q_T H^{\sigma-2/q}_r \| \}.
\end{equation}
Here we have treated $\sum_{j=1}^3 g_j$ 
in the same way as $2i \A\cdot \nabla v+|\A|^2 v+\phi v$ in the proof of Lemma~\ref{lem:smoothing}.
Converting \eqref{eq:diff2} into integral form by the use of the propagator $U(t,\tau)$ for \eqref{eq:LS},
taking the $L^\infty_T H^{s-2}$-norm and using \eqref{eq:equiv3}, we obtain
\begin{align} 
\| u_{-};\Sigma^{1,s}_T \|&\le C(R)\bigl\{\| (u_0 ,\A_0,\A_1)_{-}; X^{s,\sigma} \|
+\| u_{-};L^\infty_T L^2\|\vee \| A_{-}; L^\infty_T H^\sigma \| \nonumber \\
&\quad+\sum_{j=4}^{13}\|  g_j ;L^1_T H^{s-2}\| \bigr\}. \label{eq:auxdiffest}
\end{align}
Here we note that $K_{s}\equiv \sup_{t,\tau\in [0, T]}\| U(t,\tau); H^s\to H^s \|\le C(R)$.
We estimate the right-hand side term by term as follows:
\begin{align}
\| g_4;H^{s-2}\| &\lesssim \| \partial_t \A ;H^{\sigma-1}\| \| u_{-}; H^{s-1/2}_6 \|,  
\label{eq:g4}\\
\| g_5;H^{s-2}\| &\lesssim \| \partial_t \A ;H^{\sigma-1}\| \| \A;H^{\sigma-2/q}_r\| \| u_{-}; H^s \|, 
\label{eq:g5}\displaybreak[0] \\
\| g_6;H^{s-2}\| &\lesssim \| u;\Sigma^{1,s}_T\|^2 \| u_-; H^s \|, 
\label{eq:g6}\\
\| g_7;H^{s-2}\| &\lesssim \| \partial_t \A_{-} ;H^{\sigma-1}\| \| u'; H^{s-1/2}_6 \|, 
\label{eq:g7}\displaybreak[0] \\
\| g_8;H^{s-2}\| &\lesssim \| \A_{-};H^{\sigma-2/q-\delta}_r\| \| \partial_t u'; H^{s-1} \|, 
\label{eq:g8}\\
\| g_9;H^{s-2}\| &\lesssim \| \partial_t \A_{-} ;H^{\sigma-1}\| \| \A_+;H^{\sigma-2/q}_r\| \| u'; H^s \|,
\label{eq:g9}\displaybreak[0] \\
\| g_{10};H^{s-2}\| &\lesssim \| \A_{-};H^{\sigma-2/q}_r\| \| \partial_t \A_+ ;H^{\sigma-1}\| \| u'; H^s \|, 
\label{eq:g10} \\
\| g_{11};H^{s-2}\| &\lesssim \| \A_{-};H^{\sigma-2/q-\delta}_r \cap H^1\| 
	\| \A_+ ;H^{\sigma-2/q}_r \cap H^1\| \| \partial_t u'; H^{s-1} \|,
\label{eq:g11} \displaybreak[0] \\
\| g_{12};H^{s-2}\| &\lesssim \| u_{-} ;\Sigma^{1,s}_T\| \| u_+;\Sigma^{1,s}_T\| \| u'; H^s \|, 
\label{eq:g12}\\
\| g_{13};H^{s-2}\| &\lesssim \| u_{-} ;H^s\| \| u_+;H^s\| \| u'; \Sigma^{1,s}_T \|. \label{eq:g13}
\end{align}
We remark that we can obtain the estimates above for $s\ge 1$ and do not need the assumption $s\ge11/8$. 
In the proof of \eqref{eq:g4}-\eqref{eq:g13}, we mainly use the Leibniz rule if $s>2$, 
the H\"older and the Sobolev inequalities, Lemma~\ref{lem:phiest}, 
and the inclusions $H^{\sigma-1}\embedding L^q$ and $H^{\sigma-2/q-\delta}_r \embedding L^\infty$
together with duality argument if necessary.
For example, we can obtain \eqref{eq:g4} for $s=1,2$ by the use of the inclusion $H^{s-1/2}_6\embedding H^{s-1}_r$, 
and for $1<s<2$ by interpolation. For $s>2$, by the Leibniz rule we obtain 
\[ \| g_4;H^{s-2}\| \lesssim \| \partial_t \A ;H^{s-2}\| \| u_{-}\|_\infty
+\| \partial_t \A \|_q \| u_{-}; H^{s-2}_r \|,  \]  
and hence we obtain \eqref{eq:g4} by using the tools mentioned above.
We can analogously estimate \revised{$g_5, g_7, g_9$ and $g_{10}$}. 
We next estimate $g_8$. For $s>2$,
\begin{equation}\label{eq:auxg8}
\| g_8:H^{s-2}\| \lesssim \|\A_{-} \|_\infty \| \partial_t u';H^{s-1} \|
+\|\A_{-} ;H^{s-2}_{p_1}\| \| \nabla \partial_t u' \|_{p_2}
\end{equation}
by the Leibniz rule, where $p_1=\max\{r;3/(s-2)\}$ and $1/p_2=1/2-1/p_1$
so that $H^{\sigma-2/q-\delta}_r \embedding H^{s-2}_{p_1}\cap L^\infty$
and $H^{s-2}\embedding L^{p_2}$.
For $1\le s\le 2$, we can prove \eqref{eq:auxg8} without the second term in the right-hand side 
similarly as in the estimate of $g_4$. 
Therefore we obtain \eqref{eq:g8}.
\revised{We can analogously estimate $g_{11}$.}
The estimates for $g_6,g_{12},g_{13}$ are easy.
We should also estimate $\| u_{-};L^\infty_T L^2\|$ in \eqref{eq:auxdiffest}. 
This can be done as in the proof of Proposition~\ref{prop:localexistence}, namely the inequality
\begin{equation}
\| u_{-};L^\infty_T L^2\| \vee \| \A_{-};L^\infty_T H^{1/2}\cap L^4_T L^4\|
\le C(R) \| (u_0 ,\A_0,\A_1)_{-}; X^{0,1/2} \| \label{eq:diffL2}
\end{equation}
is obtained for sufficiently small $T$.
Applying \eqref{eq:g4}-\eqref{eq:g13} and \eqref{eq:diffL2} to \eqref{eq:auxdiffest}, using the H\"older inequality for the time variable,
and choosing $T$ sufficiently small, we obtain
\begin{align}
\| u_{-};\Sigma^{1,s}_T \|\le C(R)\bigl\{ &\| (u_0 ,\A_0,\A_1)_{-}; X^{s,\sigma} \|+
\|  \A_{-} ;\M^{1,\sigma}_T \cap L^q_T H^{\sigma-2/q}_r\|  \nonumber \\
&+\|  \A_{-} ; L^\infty_TH^1\cap L^q_T H^{\sigma-2/q-\delta}_r\| \| u';\Sigma^{1,s+1}_T\|\bigr\}.
\label{eq:auxdiffest2}
\end{align}
Next we estimate the Maxwell part.
By applying Lemma~\ref{lem:StKG} to the equation of the difference $\A_{-}$,
\begin{align*}
\|  \A_{-} ;\M^{1,\sigma}_T \cap L^q_T H^{\sigma-2/q}_r\|  
&\lesssim \| (\A_0,\A_1)_{-}; H^\sigma \oplus H^{\sigma-1} \| +\| \A_{-};L^1_T H^{\sigma-1}\| 
\nonumber \\
&\quad+\| P\J_{-};L^1_T H^{\sigma-1}\|.
\end{align*}
We have the expression
\begin{equation}\label{eq:diffJ}
P\J_{-}= 2\Im P\Bar{u}_+ \nabla u_{-} -2P(\A\Re (\Bar{u}_{+} u_-)) -2P(\A_{-} |u'|^2).
\end{equation}
Therefore a slight modification of Lemma~\ref{lem:nlmaxwell} shows 
\begin{align*}
&\| P\J_{-};L^1_T H^{\sigma-1}\| \\
&\quad \lesssim T^{1/4}\{ 
\| u_{-};L^\infty_T H^s \cap L^2_T H^{s-1/2}_6 \| \| u_+;L^\infty_T H^s \cap L^2_T H^{s-1/2}_6 \|
\langle \| \A;L^\infty_T H^\sigma \| \rangle \\
&\qquad\qquad+ \| u';L^\infty_T H^s \|^2 \| \A_{-} ; L^\infty_T H^\sigma \| \} \\
&\quad\le C(R) T^{1/4}\{ \| u_{-};L^\infty_T H^s \cap L^2_T H^{s-1/2}_6 \| 
\vee  \| \A_{-} ; L^\infty_T H^\sigma \| \}.
\end{align*}
Choosing $T$ sufficiently small,  we obtain
\begin{align*}
&\|  \A_{-} ;\M^{1,\sigma}_T \cap L^q_T H^{\sigma-2/q}_r\|   \\
&\quad \le C(R) \{ \| (\A_0,\A_1)_{-}; H^\sigma \oplus H^{\sigma-1} \| 
+T^{1/4}\| u_{-};L^\infty_TH^s \cap L^2_T H^{s-1/2}_6 \| \},
\end{align*}
which is \eqref{eq:diffest2}. 
Substituting this inequality into \eqref{eq:auxdiffest2}, we can also prove \eqref{eq:diffest1}. 
We proceed to the proof of \eqref{eq:diffest3}.
For the Schr\"odinger part, we can prove
\begin{align}
\| u_{-};\Sigma^{1,s-1}_T \|\le C(R)\bigl\{ &\| (u_0 ,\A_0,\A_1)_{-}; X^{s-1,\sigma-\delta} \|+
\|  \A_{-} ;\M^{1,\sigma-\delta}_T \cap L^q_T H^{\sigma-2/q-\delta}_r\|\bigr\}.
\label{eq:auxdiffest3}
\end{align}
If $s>2$, we can prove \eqref{eq:auxdiffest3} similarly as \eqref{eq:auxdiffest2}.
Indeed, for the estimates \eqref{eq:diffsmoothing}-\eqref{eq:g13} except \eqref{eq:g8} and \eqref{eq:g11},
we can replace $s$ with $s-1$, and  $\sigma$ with $\sigma-\delta$ respectively
since $(s-1,\sigma-\delta)\in \mathcal{R}_*$.
On the other hand, for \eqref{eq:g8} and \eqref{eq:g11} we replace $s$ with $s-1$ to obtain
\begin{align*}
\| g_8;H^{s-3}\| &\lesssim \| \A_{-};H^{\sigma-2/q-\delta}_r\| \| \partial_t u'; H^{s-2} \|, 
\\
\| g_{11};H^{s-3}\| &\lesssim \| \A_{-};H^{\sigma-2/q-\delta}_r \cap H^1\| 
	\| \A_+ ;H^{\sigma-2/q}_r \cap H^1\| \| \partial_t u'; H^{s-2} \|.
\end{align*}
In the replacements above, we do not meet the harmful factor $\| \partial_t u'; H^{s-1} \|$.
Therefore we obtain \eqref{eq:auxdiffest3}.
On the other hand, if $s\le 2$, we directly estimate $\| u_{-}; L^\infty_T H^{s-1}\|$.
To this end, we estimate the $L^1_T H^{s-1}$-norms of $g_1,g_2$ and $g_3$ in \eqref{eq:diff}
similarly as in the proof of Lemma~\ref{lem:smoothing}.
Indeed we can show
\begin{align}
\| g_1;L^1_T H^{s-1}\| &\lesssim T^{1/2} 
\| \A_{-} ; L^\infty_T H^{\sigma-\delta}\cap L^q_T H^{\sigma-2/q-\delta}_r \| 
\| u';L^\infty_T H^s \cap L^2_T H^{s-1/2}_6 \|, 
\\
\| g_2;L^1_T H^{s-1}\|&\lesssim T^{1/2} 
\| \A_{-}; L^\infty_T H^{\sigma-\delta}\cap L^q_T H^{\sigma-2/q-\delta}_r \|  \nonumber \\
&\qquad\times\| \A_+; L^\infty_T H^{\sigma-\delta}\cap L^q_T H^{\sigma-2/q-\delta}_r \| 
\| u';L^\infty_T H^s \|,  
\\
\| g_3;L^1_T H^{s-1}\| &\lesssim T \| u_{-};L^\infty_T H^{s-1}\|\| u_+;L^\infty_T H^s\|
\| u';L^\infty_T H^s\|.
\end{align}
We also have the estimate
\[
\|\partial_tu_- ; L^\infty_TH^{s-3}\|\le C(R)
\left(\|u_- ; L^\infty H^{s-1}\|+\|\A_- ; L^\infty H^1\|\right)
\]
using \eqref{eq:diff} and Lemma \ref{lem:phiest}.
These estimates prove \eqref{eq:auxdiffest3} for $1\le s \le 2$.
For the Maxwell part, we can show 
\begin{align}
&\| \A_{-};\M^{1,\sigma-\delta}_T\cap L^q_T H^{\sigma-2/q-\delta}_r \|  
\nonumber \\
&\quad\le C(R) \{ \| (\A_0,\A_1)_{-}; H^{\sigma-\delta} \oplus H^{\sigma-1-\delta} \| 
+T^{1/r} \| u_{-}; L^\infty_T H^{s-1} \cap L^2_T H^{s-3/2}_6 \| \}. \label{eq:auxdiffest4}
\end{align}
To show \eqref{eq:auxdiffest4}, we should estimate $P\J_-$ written in the form \eqref{eq:diffJ} term by term.
We first estimate $P(\Bar{u}_+\nabla u_-)$ \revised{in the case $0<\theta\equiv\sigma-1-\delta\le 1$}.
We use the inequality
$\| P(\Bar{u}_+\nabla u_-); H^\theta \| \lesssim \| \nabla u_+\|_4 \| u_-; H^\theta_4 \|$,
which is directly proved for $\theta=0,1$ and generalized for $0<\theta<1$ by interpolation.
Then we have
\[ \| P(\Bar{u}_+\nabla u_-); H^{\sigma-1-\delta} \| 
\lesssim \| u_+;H^{s-3/8}_4\| \| u_{-} ; H^{s-11/8}_4\|. \]
We remark that we have used here the assumption $s>11/8$ 
\revised{since $0<\sigma-1-\delta \le s-11/8$} (and we do not use this assumption elsewhere).
\revised{This inequality can be proved by the use of Lemma~\ref{lem:commutator} even if $\theta>1$.}
We next estimate the term $\A \Re (\Bar{u}_+ u_-)$ again by the Leibniz rule and the Sobolev inequality:
\begin{align*}
\| \A \Re (\Bar{u}_+ u_-); H^{\sigma-1-\delta} \| 
&\lesssim \| \A ; H^{\sigma-1-\delta}_{p_3} \| \| u_{+} \|_\infty \|u_- \|_{p_4}
+\| \A\|_\infty \| u_+;H^{\sigma-1-\delta}_{p_3} \| \|u_- \|_{p_4} \\
&\quad+\| \A\|_\infty \| u_+\|_\infty \|u_- ;H^{\sigma-1-\delta}\| \\
&\lesssim \| \A ; H^{\sigma-2/q-\delta}_r \| \| u_+;H^{s-1/2}_6 \| \| u_-; H^{s-1} \|.
\end{align*}
Here $1/p_3=\delta/3$ and $1/p_4=1/2-\delta/3$.
We can analogously treat the term $\A_{-} |u'|^2$.
Collecting these estimates, using \revised{the fact $H^{s-3/8}_4=(H^s, H^{s-1/2}_6)_{[3/4]}$,
and} the H\"older inequality for the time variable, and 
choosing $T$ sufficiently small, we obtain \eqref{eq:auxdiffest4}.
Substituting \eqref{eq:auxdiffest4} into  \eqref{eq:auxdiffest3}, we obtain \eqref{eq:diffest3}
\QED

\begin{proposition}\label{prop:continuous}\samepage
Let $T>0$\textup{,} $\revised{s>11/8}$\textup{,} $\sigma>1$ and let $(s,\sigma)\in\mathcal{R}$ with 
$(s+1,\sigma)\in\mathcal{R}_*$\textup{.}
Then\textup{,} the mapping defined by $(u_0,\A_0,\A_1) \mapsto (u,\A, \partial_t \A)$
is continuous as a mapping from $X^{s,\sigma}$ to $C_T X^{s,\sigma}$\textup{.}
\revised{Here\textup{,} $(u,\A)$ is the solution to
\textup{(MS-C)} with \eqref{eq:IDC} obtained in 
Propositions~\textup{\ref{prop:localexistence}-\ref{prop:regularity}. }}
\end{proposition}
\Proof
\revised{We may assume $s<4$ since the case $s\ge 4$  has already been proved in \cite{NW05}.}
Let $\eta$ be a rapidly decreasing function on $\R^3$ satisfying $\int \eta (x) dx =1$,
and let $\eta_\epsilon =\epsilon^{-3}\eta (\cdot/\epsilon)$.
We put $u_0^\epsilon =\eta_\epsilon * u_0$ and $\A_j^\epsilon =\eta_{\epsilon^{1/\delta}} * \A_j$, $j=0,1$,
and let $(u^\epsilon, \A^\epsilon)$ be a corresponding solution.
Then for $j=0,1$,
\[ \| (u_0^\epsilon,\A_0^\epsilon,\A_1^\epsilon) ;X^{s+j,\sigma}\|=O(\epsilon^{-j}), \quad
\| (u_0-u_0^\epsilon,\A_0-\A_0^\epsilon,\A_1-\A_1^\epsilon) ;X^{s-j,\sigma-j\delta}\|=o (\epsilon^j) \] 
as $\epsilon \downarrow 0$. 
We also have $\| u^\epsilon;\Sigma^{1,s+1}_T \|\lesssim \| u_0^\epsilon;H^{s+1}\|$ by Lemma~\ref{lem:Hsest}.
We use \eqref{eq:diffest1} and \eqref{eq:diffest3} for bounding the term 
$\|u-u';H^s\|$ taking $(s,\sigma)$ in Lemma \ref{lem:diffest} as 
$(s,\sigma_0)$ with $1<\sigma_0\le \sigma$, $(s,\sigma_0)\in \mathcal{R}$, 
$(s+1,\sigma_0)\in \mathcal{R}_\ast$, \revised{and $0<\sigma_0-1-\delta\le s-11/8$}.
On the other hand, we use \eqref{eq:diffest2} for bounding the term 
$\|(\A-\A^\varepsilon, \partial_t\A-\partial_t\A^\varepsilon) ; H^\sigma\oplus H^{\sigma-1}\|$ for $(s,\sigma)$.
Then we can obtain  
\begin{align*}
&\| (u-u^\epsilon,\A-\A^\epsilon,\partial_t\A-\partial_t \A^\epsilon) ;X^{s,\sigma}\| \\
&\quad \le C(R)\{ \| (u_0-u_0^\epsilon,\A_0-\A_0^\epsilon,\A_1-\A_1^\epsilon) ;X^{s,\sigma}\| \\
&\qquad+\| (u_0-u_0^\epsilon,\A_0-\A_0^\epsilon,\A_1-\A_1^\epsilon) ;X^{s-1,\sigma-\delta}\| 
\| u_0^\epsilon;H^{s+1}\| \} \\
&\quad =o(1)+o(\epsilon)O(\epsilon^{-1}) =o(1),
\end{align*}
which proves that  $(u^\epsilon, \A^\epsilon,\partial_t \A^\epsilon)$ converges to 
$(u,\A,\partial_t \A)$ in $C_T X^{s,\sigma}$.
Next we consider a sequence  $\{(u_0^n,\A_0^n,\A_1^n)\}_{n=1}^\infty$  converging to $(u_0,\A_0,\A_1)$
in $X^{s,\sigma}$.
We shall prove that the corresponding sequence of 
the solutions $\{ (u^n, \A^n,\partial_t \A^n)\}_{n=1}^\infty$ 
converges to $(u,\A,\partial_t \A)$ in $C_T X^{s,\sigma}$, which is the assertion of the proposition.
By the previous step, $(u^{n\epsilon}, \A^{n\epsilon},\partial_t \A^{n\epsilon})$ 
converges to $ (u^n, \A^n,\partial_t \A^n)$ in $C_T X^{s,\sigma}$ uniformly with respect to $n$
as $\epsilon \downarrow 0$. 
Moreover for any fixed $\epsilon$, $(u^{n\epsilon}, \A^{n\epsilon},\partial_t \A^{n\epsilon})$ converges to 
$(u^\epsilon, \A^\epsilon,\partial_t \A^\epsilon)$ as $n\to \infty$ by virtue of \eqref{eq:diffest3},
since they are sufficiently smooth solutions.
Thus we can prove the convergence of  $\{ (u_n, \A_n,\partial_t \A_n)\}_n$ by standard argument.
\QED


\section{The cases of Lorentz and the temporal gauges}
\Proofof{Theorems~\textup{\ref{thm:MS-L}} and ~\textup{\ref{thm:MS-T}}}
We prove the theorems by the gauge transform.
For any solution $(u^\LG,\phi^\LG,\A^\LG)$ to (MS-L), there exists a solution to (MS-C)
which is gauge equivalent to $(u^\LG,\phi^\LG,\A^\LG)$.
Indeed, let us put $\lambda=\Delta^{-1}\Div\A^\LG$, $\A^\CG=P\A^\LG$, 
$\phi^\CG=(-\Delta)^{-1}\rho(u^\LG)$ and $u^\CG=e^{-i\lambda}u^\LG$.
Then $(u^\LG,\phi^\LG,\A^\LG)$ and $(u^\CG,\phi^\CG,\A^\CG)$ are connected by the relation \eqref{eq:gauge}, 
and $(u^\CG,\A^\CG)$ satisfies (MS-C).
Therefore we can prove Theorem~\textup{\ref{thm:MS-L}} by Theorem~\ref{thm:MS-C}.
The assumption $\sigma\ge s-1$ is needed to ensure 
the solution to (MS-C) obtained by the gauge transform having the desired regularity.
The case of the temporal gauge can be treated analogously.
For detail, see~\cite[\S 8]{NW05}. 
\QED


 \end{document}